\documentclass[journal]{IEEEtran}
\usepackage{cite}
\usepackage{url}
\usepackage{empheq}
\usepackage{float}
\usepackage{dsfont}
\usepackage{amsthm}
\usepackage{algorithm}
\usepackage{algorithmic}

\usepackage{amsmath,amssymb,amsfonts}
\allowdisplaybreaks[4]
\usepackage{graphicx}
\usepackage{textcomp}
\usepackage{amsmath}
\usepackage{enumerate}
\usepackage{epstopdf}
\usepackage{array}
\usepackage{booktabs}
\usepackage{subfigure}
\usepackage{multirow}
\usepackage{soul, color}
\usepackage[usenames,dvipsnames]{xcolor}
\usepackage[latin1]{inputenc}

\newtheorem{proposition}{Proposition}

\newcommand{\beq}{\begin{equation}}
\newcommand{\eeq}{\end{equation}}
\newcommand{\bsq}{\begin{subequations}}
\newcommand{\esq}{\end{subequations}}

\soulregister\cite7 % Õë¶Ô\citeÃüÁî
\soulregister\citep7 % Õë¶Ô\citepÃüÁî
\soulregister\citet7 % Õë¶Ô\citetÃüÁî
\soulregister\ref7 % Õë¶Ô\refÃüÁî
\soulregister\pageref7 % Õë¶Ô\pagerefÃüÁî

\usepackage{nomencl}
\makenomenclature
\usepackage{etoolbox}
\renewcommand\nomgroup[1]{%
  \item[\bfseries
  \ifstrequal{#1}{A}{Acronyms}{%
  \ifstrequal{#1}{S}{Symbols}{%
  \ifstrequal{#1}{U}{Units}{}}}%
]}

\begin{document}

\title{Optimal Real-time Bidding Strategy For EV Aggregators in Wholesale Electricity Markets}

\author{
Shihan Huang,
Dongkun Han,
John Zhen Fu Pang,
and Yue Chen,~\IEEEmembership{Member,~IEEE}
%\thanks{This work was supported by the Chinese University of Hong Kong (CUHK) Direct Grant for Research No. 4055169. (Corresponding to Y. Chen)}
%\iffalse
\thanks{S. Huang, D. Han and Y. Chen are with the Department of Mechanical and Automation Engineering, the Chinese University of Hong Kong, HKSAR, China (e-mail: shhuang@link.cuhk.edu.hk, dkhan@mae.cuhk.edu.hk, yuechen@mae.cuhk.edu.hk) (Corresponding to Y. Chen).
}
\thanks{J. Z. F. Pang is with the Institute of High Performance Computing (IHPC), A*STAR, Singapore 138632, Singapore, (email:
john\_pang@ihpc.a-star.edu.sg).}}
\markboth{Journal of \LaTeX\ Class Files,~Vol.~XX, No.~X, Feb.~2019}%
{Shell \MakeLowercase{\textit{et al.}}: Bare Demo of IEEEtran.cls for IEEE Journals}

\maketitle

\begin{abstract}%(PARI) problem, approach, result, impact
With the rapid growth of electric vehicles (EVs), EV aggregators have been playing a increasingly vital role in power systems by not merely providing charging management but also participating in wholesale electricity markets. This work studies the optimal real-time bidding strategy for an EV aggregator. Since the charging process of EVs is time-coupled, it is necessary for EV aggregators to consider future operational conditions (e.g., future EV arrivals) when deciding the current bidding strategy. However, accurately forecasting future operational conditions is challenging under the inherent uncertainties. Hence, there demands a real-time bidding strategy based solely on the up-to-date information, which is the main goal of this work. We start by developing an online optimal EV charging management algorithm for the EV aggregator via Lyapunov optimization. Based on this, an optimal real-time bidding strategy (bidding cost curve and bounds) for the aggregator is derived. Then, an efficient yet practical algorithm is proposed to obtain the bidding strategy. It shows that with the proposed bidding strategy, the aggregator's profit is nearly offline optimal. Moreover, the wholesale electricity market clearing result aligns with the individual aggregator's optimal charging strategy given the prices. 
Case studies against several benchmarks are conducted to evaluate the performance of the proposed method.
\end{abstract}

\begin{IEEEkeywords}
EV aggregator, real-time bidding, electricity market, Lyapunov optimization, optimal charging
\end{IEEEkeywords}

\IEEEpeerreviewmaketitle

\section*{Nomenclature}
\addcontentsline{toc}{section}{Nomenclature}

% \subsection{Acronyms}
% \begin{IEEEdescription} [\IEEEusemathlabelsep \IEEEsetlabelwidth{$E_v^{\min}/E_v^{\max}$}]
% \item [DER] {Distributed energy resource}
% \item [EV] {Electric vehicle}
% \item [MPC] {Model predictive control}
% \item [OPF] {Optimal power flow}
% \item [SOC] {State of charge}
% \end{IEEEdescription}

\subsection{Parameters}
\begin{IEEEdescription}[\IEEEusemathlabelsep \IEEEsetlabelwidth{$E_v^{\min}/E_v^{\max}$}]
% \item [$A_g$] {The maximum arrival rate of charging tasks of group $g$.}
\item [$\eta$] {Charging efficiency.}
\item [$E_v^{\text{a}}, E_v^{\text{d}}$] {Energy levels of EV $v$ upon arrival and departure.}
\item [$E_v^{\min}/E_v^{\max}$] {Minimum/maximum energy level of EV $v$.}
\item [$G$] {Number of groups of EVs.}
\item [$P_v$] {Maximum charging power of EV $v$.}
\item [$P_v^{\text{end}}$] {Charging power of EV $v$ in the last time slot of charging.}
\item [$R_g$] {Maximum allowable charging delay of EVs in group $g$.}
\item [$T$] {Number of time slots.}
\end{IEEEdescription}

\subsection{Variables}
\begin{IEEEdescription}[\IEEEusemathlabelsep]
%\item [$a_g, a_v$] {Arrival rate of charging tasks of group $g$ and EV $v$.}
\item [$e_v$] {Energy level of EV $v$.}
\item [$p_v$] {Charging power of EV $v$.}
\item [$\pi$] {Electricity price.}
\item [$q_g, z_g$] {Virtual queues of group $g$.}
\item [$\mathcal{S}, \mathcal{S}_g$] {The set EVs served by an aggregator, and the set of EVs in group $g$.}
\item [$x_g$] {Charging power of group $g$.}
\end{IEEEdescription}

\section{Introduction}

As the demand for decarbonization intensifies, the prevalence of electric vehicles (EVs) have gained significant momentum in modern society. The International Energy Agency predicts that there will be 244 million EVs (including battery and plug-in electric vehicles) and more than 12 million charging points globally by 2030 \cite{iea}. EV aggregators play a significant role in the transition towards efficient and low-carbon energy systems by coordinating small-scale EV charging and interacting with the grid on their behalf. Excitingly, since 2020, based on FERC Order 2222 \cite{ferc}, EV aggregators have been allowed to participate directly in the wholesale electricity market in the U.S.

%In a real-time electricity market, EV aggregators can interact with the grid operator by real-time bidding, which facilitates the optimal energy dispatch of power grids and the integration of renewable energy sources. With the existence of fluctuating electricity price and renewable energy generation, EV aggregators optimizes EV charging patterns to maximize energy efficiency and maintain the quality of service. Meanwhile, they contribute to load management and peak shaving of grids, whereby peak electricity demand can be avoided \cite{Lyu}.
% Reorganize the literature review following this: (1) many literature studies optimal energy management of aggregators (mainly offline); (2) as aggregators are allowed to participate in electricity markets, more and more literature started to look at the market design and optimal bidding of aggregators; (3) However, the above energy management and bidding literature relied on predictions, which is impractical. There are online algorithms (MPC, Lyapunov optimization, learning, ...). (4) But these work focused on energy management problem and seldom consider the optimal bidding problem. The bidding problem is more complicated than the energy management in that xxxx

% Energy management or pricing scheme without bidding: 
Optimal energy management of EV charging by aggregators, utilizing a wide range of tools, have been studied in the literature. A distributed optimization framework for EV aggregators using alternating directions method of multipliers was developed in \cite{Rivera}, with the goals of demand-side management and cost minimization. A comprehensive review of distributed algorithms applied to EV charging was provided by \cite{Nimalsiri}. 
On the other hand, game theoretical methods have also been used to study the interplay between aggregators and EV drivers, such as in \cite{Alsabbagh} who utilized a non-cooperative Stackelberg game, and \cite{Gupta} who studied pricing of EV charging under the aggregator setting.
The studies above focus on how the aggregator set charging prices for EVs.
%Moreover, most of them derived the optimal strategy of aggregator from offline optimization problems.

%As aggregators are allowed to participate in electricity markets, more and more literature started to look at the market design and optimal bidding of aggregators in wholesale markets. 
Driven by FERC order 2222, increasing attention has been given to market design and optimal bidding of aggregators in the electricity market. 
Reference \cite{Duan} examined the bidding strategies for an aggregator in both energy and reserve markets.
% A mixed integer quadratic programming problem was formulated in \cite{Duan} for an aggregator with fast charge stations and energy storage in energy and reserve markets. 
% {\color{cyan} Although relaxed to a mixed integer linear programming problem, an iterative solution was still necessary due to computation complexity.}
The worst-case EV availability in terms of battery draining and energy exchange was considered in \cite{Porras}. Step bidding curve was exploited in \cite{Wang} by the conditional value-at-risk method.
A market settlement framework for distributed energy resource (DER) aggregators in a three-phase unbalanced distribution network was developed in \cite{Mousavi}.
An iterative bidding framework was proposed in \cite{Hou} by game theoretic analysis and parallel machine scheduling model. The bidding strategy of a virtual power plant with a large number of EVs was derived \cite{YangD}. However, the above literature focus on the deterministic setting, overlooking uncertainties, e.g., random renewable supply. Additionally, the optimization problems established in most of the above are non-convex, and without closed-form solutions; as a result, iterative methods are necessary, and these may either not converge or have a high computational burden.
Apart from the uncertainties present in the electricity market e.g., supply, demand, and electricity prices, the EV charging network brings about new sources of uncertainty, such as in the arrival and departure of EVs. To address these, a two-stage stochastic optimal bidding model for EV aggregators in day-ahead and real-time energy and regulation markets was proposed in \cite{Vagropoulos2013}, while uncertain market prices and driving requirements were similarly addressed in \cite{BARINGO2017362}. Further, to consider these uncertainties under a game-theoretic setting, a stochastic generalized Nash game was formulated in \cite{moghaddam2018network}. While multiple works, e.g., \cite{Vagropoulos2013, Lyu, Han, Zhou, Song}, employed stochastic programming to tackle uncertainty, they require predictions, i.e., assuming some information about the future. In practice, some uncertain factors of EVs are hard to predict at all, e.g., future incoming EVs. Therefore, an online algorithm without the need for predictions is desired.

%based on which a centralized real-time framework of EV charging management was proposed \cite{Vagropoulos2016}.

% Therefore, online algorithms are essential in the research of optimal bidding of EV aggregators, enabling the response to the dynamic and changing conditions of the EV market and real-time decision making. 
% However, in most of the existing research, electricity price was either considered as an input parameter %, e.g. \cite{Mousavi, Rivera}, 
% or forecast based on historical data. %, e.g. \cite{YangD, Liu}. 
% \cite{YangD}
%The mechanism of electricity market, especially how the price is determined in bidding, has not been well addressed. The optimal bidding strategy derived by the model also depended on the forecast electricity prices in the day-ahead and real-time markets, although the reliability of the prediction is open to question. 
% It is common in the literature to utilize the prediction of uncertainties in the system, 
Online algorithms have attracted great attention in recent years due to the  need to address uncertainties. Model predictive control (MPC) has been widely used in the online scheduling \cite{Liu} and bidding \cite{Nakano} of EV aggregators. However, the MPC method still relies partly on predictions. Another class of online algorithms is the learning-based algorithms \cite{Gao,najafi2019reinforcement,tao2021deep}. However, to achieve a satisfactory performance by these learning-based algorithms, a large amount of data is typically necessary for training and testing, which increases the cost of deployment. Additionally, the EV aggregators may not always have sufficient data, especially newly built ones.
% In reference \cite{Song}, day-ahead price was directly input into the model, while real-time price was forecast by autoregressive integrated moving average and Markov process. Historical market data was used in \cite{Lyu} based on the assumption that the prices could be accurately forecast.
% EV batteries were grouped in \cite{Han} according to their departure time and stages, and their driving routes were predicted, whereby a multi-stage stochastic optimization problem was formulated. The optimal bidding strategy of EV aggregator was modelled in \cite{YangH} by an expectation minimization problem with conditional value-at-risk constraints considering the uncertainty of electricity price. A risk-constrained two-stage stochastic optimization model was proposed in \cite{Zhou} by the introduction of the dispatchable region of large-scale EV fleet. However, the accuracy of the dispatchable region depends on the uncertainty of EVs and distributed energy resources (DERs), which requires reliable prediction of EV behaviors and renewable generations. A bilevel optimization problem was formulated in \cite{Porras} by robust optimization and subsequently transformed into single-level mixed-integer linear programming, whereas limited to the operations in the day-ahead market. Hybrid stochastic-robust optimization, which is a combination of stochastic programming and robust optimization, has been applied in \cite{MOGHADDAM2018478, BARINGO2017362}.
Moreover, works that utilize learning-based algorithms cannot provide an explicit bidding strategy with theoretical guarantee and the EV aggregators may lack trust in the results.
%Other innovative approaches have been also adopted in the literature recently, such as learning-based algorithm \cite{Gao}, fuzzy logic \cite{Li}, etc. 
%A stochastic game with oligopoly players wind power plants and EV aggregators is formulated by multi-agent reinforcement learning in \cite{Gao}. 
% However, only bidding decisions in the day-ahead market were investigated and those in the real-time market were not addressed. A dynamic pricing model was developed in \cite{MoghaddamZ} and the optimal strategies of charging stations were obtained by ant colony optimization.
% Charging dispatch was modelled by fuzzy decision making in \cite{Li}, thereby a decentralized proximal best response algorithm was developed. 
%Despite the advantages of these novel techniques, they have some noteworthy limitations. For example, they can be computationally expensive and time-consuming for real-time decision making in energy systems, especially on large-scale model. The results of these algorithms may be sensitive to the choice of hyperparameters or fuzzy sets. Large data sets are also necessary for learning-based algorithms to obtain satisfying performance, which increases the cost of deployment.
Lyapunov optimization is another class of online methods that can offer near-optimal strategies without the need for predictions in the system. % By defining an appropriate Lyapunov function, the stability of a dynamic system can be analyzed without explicitly modelling its dynamics, thereby a optimal control policy is derived by drift-plus-penalty algorithm.
Unlike prediction-based techniques such as MPC, Lyapunov optimization does not rely on the prediction of uncertainties, %such as renewable generations, arrivals of EV charging demands, energy level of EVs, etc.
thus particularly suitable for energy system integrated with DERs and EVs, which possess highly random behaviors.
The existing applications of Lyapunov optimization concentrated on the optimal energy management for EV charging stations \cite{zhou2017optimal,Sun,yan2022distributed}, energy storage \cite{guo2021real}, and data centers \cite{hou2020decentralized}. However, the optimal bidding problem, which is much more complex, has not been considered. This complexity mainly lies in the fact that to derive the optimal bidding strategy, one needs to consider not only the EV charging requirements but also the potential market clearing under the bidding.

In this paper, we bridge the research gaps by developing an optimal bidding strategy for EV aggregators in a wholesale electricity market. Our main contributions are two-fold:

1)	\emph{Real-time EV aggregator bidding strategy.} To enable participation of EV aggregators in real-time electricity markets, a prediction-free bidding strategy is developed. First, the online optimal EV charging scheduling strategy by an aggregator, as a function of the real-time electricity price, is established. Based on this, the real-time bidding strategy (bidding cost curve and bounds) for an EV aggregator is derived. We prove that with the proposed bidding strategy, the EV aggregator can obtain a profit close to the offline benchmark. Moreover, electricity market clearing outcomes align with aggregators' individual interests. To the best of our knowledge, this paper is one of the first studies that provide an explicit bidding strategy for an EV aggregator with provable theoretical properties.
%1)	\emph{A prediction-free optimal bidding framework.} An optimal bidding framework of EV aggregators is proposed based on Lyapunov optimization, which does not rely on the in-advance forecast of uncertainties in the system, e.g. renewable generations, EV arrivals, etc. An optimal bidding function of aggregators is derived from the framework, which is reported to the market operator for market clearing. {\color{red} To the best of our knowledge, there are few studies covering explicit bidding functions for EV aggregators.}
% {\color{blue}The framework is designed based on the principle of privacy protection agents in the market, in which only necessary information for bidding is revealed.??}

%2) \emph{An online algorithm for real-time bidding event.} An online algorithm for the optimal bidding problem is developed for the case of quadratic bidding functions. Using the stepwise integration of the bidding function, the bidding problem can be solved programmatically even though it is difficult to derive the closed-form solution of the problem. A linearization approach is also adopted to facilitate the solution.
2) \emph{Practical algorithm to generate the bidding strategy.} The optimal bidding cost curve derived theoretically turns out to be a piecewise quadratic function. To facilitate implementations, a practical algorithm based on stepwise integration is proposed to generate the cost curve programmatically. A computationally efficient ``aggregator bidding-market clearing'' framework is established with the help of linearization techniques. Case studies demonstrate the effectiveness of the proposed algorithm and framework.

The rest of this paper is organized as follows. The online EV charging scheduling algorithm is developed in Section \ref{sec:model}. In Section \ref{sec-III}, the optimal real-time bidding strategy is derived theoretically and a practical algorithm for generating the bidding strategy is established. Case studies are conducted in Section \ref{sec:result} where the performance of the proposed method is also evaluated. Finally, conclusions are drawn in Section \ref{sec:conclu}.

\section{Online EV and Electricity Market Models} \label{sec:model}

In this section, we first formulate the EV charging scheduling problem for an EV aggregator and propose an online algorithm to solve the problem using Lyapunov optimization. Next, we introduce the electricity market as well as how aggregators participate in it.

%Based on this, the optimal bidding strategy of the aggregator in a wholesale electricity market is derived.

%the establishment of the optimal bidding framework of EV aggregators is introduced. First of all, the operation of an EV aggregator is modelled and an online cost minimization problem is formulated by Lyapunov optimization. Next, the optimal strategy of the aggregator is derived from the solution to this problem. Finally, the optimal bidding problem of aggregators is formulated.

\subsection{EV Charging Scheduling Model}

% \begin{figure}[t]
%   \centering
%   \includegraphics[width=0.45\textwidth]{Figs/Illustration1.pdf}\\
%   \caption{Illustration of the workload and energy flows among DCs.}\label{fig:sysConf}
% \end{figure}
We consider an EV aggregator that provides a set of EVs ($v \in \mathcal{S}$) with charging services. The time horizon studied is divided into $T$ time slots. Each EV $v$ arrives at the charging station at time $T_v^{\text{a}}$ and departs at $T_v^{\text{d}}$. Its initial battery energy level upon arrival is $E_v^{\text{a}}$ and its minimum required energy level of departure is $E_v^{\text{d}}$. Denote the electricity price by $\pi(t)$ and the charging power of EV $v$ by $p_v(t)$. Then, the aggregator solves the following optimization model to determine the optimal EV charging schedules with minimum cost \eqref{eq:EV-1}.
\bsq
\label{eq:EV}
\begin{align}
    \min_{p_v(t),\forall t,\forall v}~ & \sum \nolimits_{t=1}^T \left(\pi(t) \sum \nolimits_{v \in \mathcal{S}} p_v(t)\right), \label{eq:EV-1}\\
    \mbox{s.t.}~ &  0 \le p_v(t) \le P_v,\forall v,\forall t, \label{eq:EV-2}\\
    ~ & e_{v,T_v^{{\text{a}}}}= E_v^{\text{a}},~e_{v,T_v^{\text{d}}}\geq E_v^{\text{d}},\forall v,\label{eq:EV-3}\\
    ~  & e_{v,t+1}=e_{v,t}+\eta_c {p}_{v,t}^{c}\Delta t,\forall v,\forall t \ne T, \label{eq:EV-5}\\
  ~ & E^{\text{min}}_v\leq e_{v,t}\leq E^{\text{max}}_v,\forall v,\forall t, \label{eq:EV-4}
\end{align}
\esq
%\jpang{why do we need equation 1d? maybe the max to add that no need to charge beyond say 80\% or 100\%? Wouldn't the minimum be taken care by the arriving SoC already?}
where $P_v$ is the maximum charging power limited by \eqref{eq:EV-2}; $e_{v,t}$ is the battery energy level of EV $v$ with lower and upper bounds $E_v^{\text{min}}$ and $E_v^{\text{max}}$ in \eqref{eq:EV-4}; \eqref{eq:EV-3} and \eqref{eq:EV-5} give the battery energy level requirements and dynamics.

However, full information on the uncertainties (electricity prices and EV arrivals) over the whole time horizon is necessary to solve problem \eqref{eq:EV}, which is not a realistic assumption in practice. To address this issue, we provide a reformulation of \eqref{eq:EV} based on which an online algorithm can be derived. First, we divide the EVs served by the aggregator into $G$ groups according to their allowable charging delay, or equivalent, the duration of time they are parked at the charging station. Let $x_g(t)$ be the aggregate charging power of all EVs in the $g$-th group $\mathcal{S}_g \in \mathcal{S}$. Denote the arrival rate of EV charging tasks by $a_g(t)$, which satisfies:
\begin{equation}
    a_g(t) = \sum \nolimits_{v \in \mathcal{S}_g} a_v(t), \forall t
\end{equation}

The arrival rate $a_v(t)$ is determined by the EV charging need and the charging power limit. For EV $v$, $a_v(t)$ can be calculated as follows:
\begin{gather}\label{eq:agvt-lb}
    a_v(t) = \left\{\begin{array}{ll}
    P_v,     & \text{if } T_v^{\text{a}} \leq t < T_v^{\text{a}} + T_v^{\min}, \\
    P_v^{\text{end}}, & \text{if } t = T_v^{\text{a}} + T_v^{\min}, \\
    0,           & \text{otherwise},
    \end{array}\right.
\end{gather}
% % The charging demands of EVs are modelled by the method proposed in \cite{Yan}. 
% We assume that the demands of all EVs need to be satisfied as soon as possible, that is, every EV is charged at its maximum charging power until the charging is finished. Suppose the EV $v$ arrives at and leaves the charging station at the time slots $T_v^{\text{a}}$ and $T_v^{\text{d}}$ respectively, then the charging demand of the EV $v$ in the time slot $t$, denoted by $a_v(t)$, is determined as follows:
% \begin{gather}\label{eq:agvt-lb}
%     a_v(t) = \left\{\begin{array}{ll}
%     P_v,     & \text{if } T_v^{\text{a}} \leq t < T_v^{\text{a}} + T_v^{\min}, \\
%     P_v^{\text{end}}, & \text{if } t = T_v^{\text{a}} + T_v^{\min}, \\
%     0,           & \text{otherwise},
%     \end{array}\right.
% \end{gather}
where $T_v^{\min}$ is the charging time of the EV $v$ using the maximum charging power $P_v$, given by $T_v^{\min} = \left[ {E_v^{\text{c}}}/{P_v \eta} \right]$
with $[x]$ representing the least integer larger than $x$;
$P_v^{\text{end}}$ is the charging power of EV $v$ in the last time slot before it is fully charged, given by 
\begin{equation}
    P_v^{\text{end}} = (E_v^{\text{d}} - E_v^{\text{a}})/{\eta} - (T_v^{\min} - 1) P_v,
\end{equation}
where $\eta$ is the charging efficiency.

Based on this, problem \eqref{eq:EV} can be reformulated as
\bsq \label{eq:flexibility-p2}
\begin{align}
    \textbf{P1: } \min_{x_g(t),\forall g,\forall t} ~& \lim_{T\rightarrow\infty} \frac{1}{T} \sum_{t=1}^{T}\mathbb{E} \left[f(t) \right], \label{eq:flexibility-p2-obj}\\
    \text{s.t. } \quad & \lim_{T \to \infty} \frac{1}{T} \sum_{t=1}^T \mathbb{E} [a_g(t)- x_g(t)] \le 0 , \forall g, \label{ineq:ag_xg_lim} \\
    & 0 \leq x_g(t) \leq X_g(t), ~ \forall g, ~ \forall t. \label{ineq:xgub}
\end{align}
\esq
where $f(t) = \pi(t) \sum_{g = 1}^G x_g(t)$. The objective function \eqref{eq:flexibility-p2-obj} minimizes the long-term time-average cost of the aggregator. The constraint \eqref{ineq:ag_xg_lim} ensures that all charging demands are satisfied in the long run, and the charging power of each group is bounded by the constraint \eqref{ineq:xgub}. In particular, $X_g(t)$ is the maximum charging power of group $g$ and $X_g(t) = \sum_{v \in \mathcal{S}_g} X_v(t)$; and $X_v(t)$ is the maximum charging power of EV $v$ given by 

% We divide the EVs served by the aggregator into $G$ groups according to their allowable charging delay, or equivalent, the duration of time they are parked at the charging station.
% %The set of all EVs using the service is denoted by $\mathcal{S}$. 
% For the $g$-th group $\mathcal{S}_g \in \mathcal{S},~ g = 1, 2, \dots, G$, the arrival rate of charging demands $a_g(t)$ is calculated as the sum of the arrival rates of all EVs in $\mathcal{S}_g$, i.e.,
% \begin{equation}
%     a_g(t) = \sum_{v \in \mathcal{S}_g} a_v(t), \forall t
% \end{equation}

% Denote the upper bounds of the charging power of the group $g$ and the EV $v$ in the time slot $t$ by $X_g(t)$ and $X_v(t)$, respectively. Then, $X_g(t) = \sum_{v \in \mathcal{S}_g} X_v(t)$. To guarantee that the charging power of the group $g$ $x_g(t)$ can be completely disaggregated, $X_v(t)$ is determined as follows:
\begin{equation} \label{eq:xv}
    X_v(t) = \begin{cases}
    \multirow{2}{*}{$P_v,$}
    & \text{if } T_v^{\text{a}} \leq t < T_v^{\text{d}} \text{ and}\\
    & e_v(t) \le E_v^{\text{d}} - \eta P_v \Delta t, \\
    \multirow{2}{*}{$\dfrac{E_v^{\text{d}} - e_v(t)}{\eta \Delta t},$}
    & \text{if } T_v^{\text{a}} \leq t < T_v^{\text{d}} \text{ and}\\
    & e_v(t) > E_v^{\text{d}} - \eta P_v \Delta t, \\
    0, & \text{otherwise},
    \end{cases}
\end{equation}
where $\Delta t$ is the duration of a single time slot. \eqref{eq:xv} is consistent with the calculation of arrival rate of charging tasks in \eqref{eq:agvt-lb}.

%$X_v(t)$ defined by \eqref{eq:xv} and the corresponding $X_g(t)$ reflect the maximum charging power that can be allocated to the EV $v$ and the group $g$ in the time slot $t$ respectively.
%At the end of charging, the power of the EV $v$ must be limited so that its energy level will not exceed $E_v^{\text{d}}$. 
%Given the electricity price $\pi(t)$ in the time slot $t$, $t = 1, 2, \dots, T$, we can formulate the optimal EV charging problem as follows:

We introduce virtual queues $q_g(t)$ to collect the charging tasks in group $g$, which is defined by
% \eqref{ineq:ag_xg_lim} is a long-term constraint embedded with conditions over the time horizon, thus it cannot be directly solved in an online manner.
% To decompose the constraint \eqref{ineq:ag_xg_lim} into a single time slot, we establish the queues of charging demands, the backlogs of which need to be controlled. The charging power of all EVs in the group $\mathcal{S}_g$ is denoted by $x_g(t)$, then the queue of charging demands of $\mathcal{S}_g$ is defined by
\begin{equation} \label{eq:Qgt}
    q_g(t+1) = \max \{q_g(t)-x_g(t),0\} + a_g(t),
\end{equation}

With the virtual queues, constraint \eqref{ineq:ag_xg_lim} can be replaced by the mean-rate-stable constraint as follows:
\begin{equation} \label{eq:mrs_q}
    \lim_{T \to \infty} {\mathbb{E}[q_g(T)]}/{T} = 0, ~ g = 1, 2, \dots, G.
\end{equation}
Specifically, the following inequality can be derived from \eqref{eq:Qgt}:
\begin{align}\label{eq-1}
     q_g(t+1) - q_g(t) \ge a_g(t) - x_g(t), ~ \forall t.
\end{align}
Summing \eqref{eq-1} up over all time slots and dividing both sides by $T$, we have \eqref{ineq:ag_xg_lim} is met, i.e.,
\begin{align}
    \lim_{T\rightarrow\infty} \frac{1}{T} \sum_{t=1}^T \mathbb{E}[a_g(t) - x_g(t)] \le \lim_{T\rightarrow\infty} \frac{1}{T} \mathbb{E}[q_g(T+1)] = 0. \nonumber
\end{align}

Furthermore, considering that \eqref{ineq:ag_xg_lim} (or \eqref{eq:mrs_q}) only ensures that charging demands are met in a time-average sense, we introduce an additional virtual queue $z_g(t)$ to bound the charging delays:
\begin{equation} \label{eq:zgt}
    z_g(t+1) = \max \left\{z_g(t) + \frac{\alpha_g}{R_g} \mathbb{I}_g(t) - x_g(t), 0 \right\}, \forall g, \forall t,
\end{equation}
where $R_g$ is the duration of the EVs parked at the charging station in the group $\mathcal{S}_g$, i.e., $R_g = T_v^{\text{d}} - T_v^{\text{a}}$. $\alpha_g$ is a positive number that controls the behaviors of virtual queues, $\mathbb{I}_g(t)$ is an indicator function of $q_g(t)$, given by
\begin{align}
    \mathbb{I}_g(t) = \left\{\begin{array}{ll}
      0, & \text{if}~ q_g(t) = 0,\\
      1, & \text{if}~ q_g(t) > 0.
    \end{array} \right.
\end{align}
Obviously, any unserved charging demand remaining from previous time slot $t-1$ will increase $z_g(t)$ of the current time $t$. By requiring $z_g(t)$ to be mean-rate-stable as follows, we can control the charging delay.
\begin{equation} \label{eq:mrs_z}
    \lim_{T \to \infty} \frac{\mathbb{E}[z_g(T)]}{T} = 0, ~ g = 1, 2, \dots, G.
\end{equation}

%However, after we replace \eqref{ineq:ag_xg_lim} by \eqref{eq:Qgt} and \eqref{eq:mrs_q}, \textbf{P1} is still unsolvable because \eqref{eq:Qgt} is a time-coupling constraint. Hence, \textbf{P1} needs to be further reformulated.

\subsection{Online Charging Scheduling Using Lyapunov Optimization} \label{subsec:lya}
In the following, we derive an online algorithm for the EV charging scheduling problem based on Lyapunov optimization. 
% To eliminate the time-coupling constraint \eqref{eq:Qgt}, Lyapunov optimization framework is established and applied to \textbf{P1}. Compared to traditional online optimization methods, Lyapunov optimization can be implemented without the prediction of the underlying uncertainties and inherently guarantees the stability of the system. Consequently, the constraints related to queues of charging demands can be excluded from \eqref{eq:Qgt}.
%To guarantee \eqref{eq:mrs_q} and \eqref{eq:mrs_z}, we introduce Lyapunov function to apply Lyapunov optimization technique based on $q_g(t)$ and $z_g(t)$. 
\subsubsection{Procedures}
Denote the concatenated vector of queues by 
\begin{equation}
    \boldsymbol{\Theta} (t) = (q_1(t),\dots,q_G(t), z_1(t), \dots, z_G(t)).
\end{equation}
Define the Lyapunov function as
\begin{equation} \label{eq:LyaFun}
L(\boldsymbol{\Theta} (t)) \triangleq \frac{1}{2} \sum \nolimits_{g = 1}^G q_g^2(t) + \frac{1}{2} \sum \nolimits_{g = 1}^G z_g^2(t).
\end{equation}
$L(\boldsymbol{\Theta} (t))$ is a measure of the backlogs of all queues. Therefore, the increase in the backlogs of queues can be evaluated by the increment in the Lyapunov function, namely the Lyapunov drift. The Lyapunov drift from the time slot $t$ to $t + 1$ is calculated as follows:
\begin{equation} \label{eq:LyaDrift}
\begin{aligned}
\Delta(\boldsymbol{\Theta} (t)) 
&= L(\boldsymbol{\Theta}(t+1)) - L(\boldsymbol{\Theta} (t)) \\ 
&= \frac{1}{2} \sum_{g = 1}^G \left(q_g^2(t+1) - q_g^2(t)\right) \\
& \quad + \frac{1}{2} \sum_{g = 1}^G \left(z_g^2(t + 1) - z_g^2(t)\right).
\end{aligned}
\end{equation}
%where the second equality can be derived by the definitions of queues \eqref{eq:Qgt} and \eqref{eq:zgt}. 

Subsequently, the mean-rate-stable constraints \eqref{eq:mrs_q} and \eqref{eq:mrs_z} can be integrated into \textbf{P1} by modifying the objective function into a drift-plus-penalty term:
\begin{equation} \label{obj:driftPlusP}
\mathbb{E} [\Delta (\boldsymbol{\Theta} (t)) + V f(t) | \boldsymbol{\Theta} (t)],
\end{equation}
where $V$ is a positive constant to be determined. \eqref{obj:driftPlusP} strikes a balance between (virtual) queue stability and cost minimization of EV operation.

However, even after the substitution of the objective in \textbf{P1} by \eqref{obj:driftPlusP}, the problem is still time-coupled due to the term $\Delta(\boldsymbol{\Theta} (t))$. To address this, we further relax \eqref{obj:driftPlusP} into its upper bound. Particularly, for $q_g(t)$ we have
\begin{equation} \label{ineq:Q2-Ub0}
\begin{aligned}
q_g^2(t+1) &= (\max \{q_g(t)- x_g(t),0\} + a_g(t))^2 \\
            &\leq q_g^2(t) + A_g^2 + x_g^2(t) \\
            & \quad + 2q_g(t)(a_g(t)- x_g(t)),
\end{aligned}
\end{equation}
where $A_g$ is the maximum arrival of charging demands of all time slots, i.e., $A_g = \max \{a_g(t), \forall t = 1, 2, \dots, T \}$. The relaxation in \eqref{ineq:Q2-Ub0} can be applied to $z_g(t)$ similarly.
Denote the maximum backlogs of $q_g(t)$ and $z_g(t)$ of all time slots by $Q_g$ and $Z_g$ respectively, then the increments in $q_g^2(t)$ and $z_g^2(t)$ can be bounded:
\bsq
\begin{align}
    & \frac{1}{2} \left(q_g^2(t+1) - q_g^2(t)\right) \nonumber \\
    \le ~& \frac{1}{2} \left(A_g^2 + {x_g^2(t)} \right) + Q_g A_g - q_g(t) x_g(t), \label{ineq:Q2-Ub} \\
    & \frac{1}{2} \left(z_g^2(t + 1) - z_g^2(t) \right) \nonumber \\
    \le ~& \frac{1}{2} \max \left\{ \left( \frac{\alpha_g}{R_g} \right)^2, X_g^2(t) \right\} + \frac{\alpha_g Z_g}{R_g} - z_g(t) x_g(t). \label{eq:z2-lb}
\end{align}
\esq
%Note that the quadratic term $x_g^2(t)$ has been retained in \eqref{ineq:Q2-Ub0}, whereas it is usually replaced by its upper bound in the literature. We will compare the two ways of relaxation later in case studies.

Afterwards, take \eqref{ineq:Q2-Ub} and \eqref{eq:z2-lb} into the drift-plus-penalty term in \eqref{obj:driftPlusP}, yielding
\bsq
\begin{equation} \label{ineq:dpp}
    \begin{aligned}
        & \Delta(\boldsymbol{\Theta} (t)) + Vf(t) \\
        \le ~& M + V f(t) - \sum_{g = 1}^G (q_g(t) + z_g(t)) x_g(t) 
         + {\frac{1}{2} \sum_{g = 1}^G x_g^2(t)},
    \end{aligned}
\end{equation}
where $M$ is a constant given by
\begin{align}
    M = ~ & \sum_{g = 1}^G \left( \frac{A_g^2}{2} + Q_g A_g + \frac{\alpha_g Z_g}{R_g} \right) \nonumber \\
    & + \frac{1}{2} \sum_{g = 1}^G \max \left\{\left(\frac{\alpha_g}{R_g}\right)^2, \overline{X}_g^2 \right\}, \\
    \overline{X}_g = & \max \{X_g(t) ~|~ t = 1, 2, \dots, T \}, ~ \forall g.
\end{align}
\esq

Finally, we obtain the following online optimization problem by reorganizing the expression in \eqref{ineq:dpp} and ignoring the constant terms:
\bsq \label{P2}
\begin{align}
  \textbf{P2: } \min_{x_g(t),\forall g} ~ & \sum_{g = 1}^G (V\pi(t)-q_g(t) - z_g(t)) x_g(t) \nonumber \\
  & + {\sum_{g = 1}^G \frac{1}{2} x_g^2(t)}, \label{obj:EV} \\
  \hbox{s.t.}~ & 0\le x_g(t) \le X_g(t), ~ \forall g. \label{ineq:xg}
\end{align}
\esq

Compared to \textbf{P1}, the online problem \textbf{P2} is only related to the variables and parameters in time slot $t$ and can be solved in an online manner. % In each time slot, $q_g(t)$ and $z_g(t)$ are first updated based on \eqref{eq:Qgt} and \eqref{eq:zgt} before solving \textbf{P2}.

\subsubsection{Properties}
The transformation from \textbf{P1} to \textbf{P2} for the purpose of online feasibility come at the cost of sub-optimality since the objective function of \textbf{P2} is derived by drift-plus-penalty method \eqref{obj:driftPlusP}. Therefore, it is of great importance to find a bound of the optimality gap between the optimal solutions of \textbf{P1} and \textbf{P2}.
Denote the values of $f(t)$ optimized by \textbf{P1} and \textbf{P2} by $f^*(t)$ and $\hat{f}(t)$ respectively, $t = 1, 2, \dots, T$, then the optimality gap between \textbf{P1} and \textbf{P2} is bounded by the following proposition:
\begin{proposition}[Optimality] \label{prp:gap}
The optimality gap between \textbf{P1} and \textbf{P2} is bounded above, i.e.,:
\bsq
\begin{equation} \label{ineq:gap}
    \lim_{T \to \infty} \frac{1}{T} \sum_{t=1}^T \mathbb{E} \left[ f^*(t) - \hat{f}(t) \right] \le \frac{\tilde M}{V},
\end{equation}
where
\begin{align}
    \tilde M & = M + \frac{1}{2} \sum_{g = 1}^G (\overline{X}_g)^2.
\end{align}
\esq
\end{proposition}
The proof of Proposition \ref{prp:gap} can be found in Appendix \ref{appd:gap}. \eqref{ineq:gap} shows that the optimality gap decreases as $V$ increases. However, the mean-rate stability of queues \eqref{eq:mrs_q} and \eqref{eq:mrs_z} will be undermined if $V$ is too large. The resulting arger backlogs of queues are caused by the gap between charging demands and the capability of charging service provided by aggregators, which may lead to large charging delays in the meantime.
Specifically, we can show that the maximum charging delay of EVs is bounded as follows:
\begin{proposition}[Feasibility] \label{prp:delay}
The maximum charging delay of the group $g$, denoted by $D_g$, is bounded above, i.e.,:
\begin{equation} \label{ineq:delay}
    D_g \le \frac{R_g (Q_g + Z_g)}{\alpha_g}
\end{equation}
\end{proposition}
In other words, if an EV of the group $g$ arrives at the time slot $t_0$, then its charging demand will be satisfied before the time slot $t_0 + D_g$. The proof of Proposition \ref{prp:delay} can be found in Appendix \ref{appd:delay}. \eqref{ineq:delay} shows that the maximum charging delay is proportional to the sum of queues, which is consistent with intuition. 
Proposition \ref{prp:gap} and Proposition \ref{prp:delay} show that the drift-plus-penalty method strikes a balance between optimality and feasibility, which is achieved by controlling the value of $V$ in practice. Therefore, the choice of $V$ is critical to the performance of the algorithm.

\subsection{Electricity Market Model}
The above provides a mechanism for the aggregator to optimally schedule EV charging in an online manner. However, the electricity prices $\pi(t),\forall t$ in the objective of \textbf{P2} is not a given constant but determined by the electricity market with EV aggregators' bids. Specifically, each aggregator bids its utility function to the operator, based on which the electricity market is cleared by solving an optimal power flow (OPF) problem. Suppose there are $I$ generators, $J$ fixed loads and $K$ aggregators, denote the utility function and the price of the aggregator $k$ in the time slot $t$ by $u_k^{(t)}(x_k)$, the market clearing problem is as follows:
\bsq
\begin{align}
    \textbf{P3: } \min_{p_i, x_k}~ & \sum_{i = 1}^I f_i^{\text{Gen}}(p_i(t)) - \sum_{k = 1}^K u_k^{(t)}(x_k(t)) \label{obj:OPF} \\
    \mbox{s.t.}~ & \sum_{i = 1}^I p_i(t) = \sum_{j = 1}^J d_j(t) + \sum_{k = 1}^K x_k(t), ~ \forall t, \label{eq:ErgBal} \\
    & 0 \leq x_k(t) \leq X_k(t), ~ \forall k, ~ \forall t, \\
    & \underline{P}_i \le p_i(t) \le \overline{P}_i, ~ \forall i, ~ \forall t, \label{ineq:GenPwr} \\
    & -F_l \le p_l^{\text{line}} (t) \le F_l, \forall l, ~ \forall t,\label{ineq:LinePwr}
\end{align}
where $p_i(t)$ and $f_i^{\text{Gen}}(p_i(t))$ denote the power and cost of generator $i$, respectively; and $f_i^{\text{Gen}}(p_i)$ is a convex function of $p_i$. $d_j(t)$ represents the fixed demand in the network, $X_k$ is the maximum charging power of the aggregator $k$, $\underline{P}_i$ and $\overline{P}_i$ are the minimum and maximum power of the generator $i$ respectively, $F_l$ represents the maximum line power of the branch $l$.
$p_l^{\text{line}} (t)$ is the line power of the branch $l$ given by
\begin{equation} \label{eq:p_line}
    p_l^{\text{line}} (t) = \sum_{i = 1}^I \delta_{il} p_i(t) - \sum_{j = 1}^J \delta_{jl} d_j(t) - \sum_{k = 1}^K \delta_{kl} x_k(t), \forall l,
\end{equation}
\esq
where $\delta_{il}$ is the power transfer distribution factor of the node $i$ to the line $l$. 

%In this paper, we assume that $f_i^{\text{Gen}}(p_i)$  convex function of $p_i$ for each generator $i$ in the whole time horizon. Hence, \textbf{P3} is convex according to Proposition \ref{prp:con} and it can be solved efficiently.

Denote the dual variables associated with \eqref{eq:ErgBal}, the lower and upper bounds in \eqref{ineq:LinePwr} by $\lambda (t)$, $\underline{\chi}_l (t)$ and $\overline{\chi}_l (t)$ respectively, then the locational marginal price (LMP) for the aggregator $k$ in the time slot $t$ is defined as:
\begin{equation} \label{eq:price}
    \pi_k (t) = \lambda (t) + \sum_l \left( \delta_{kl} \underline{\chi}_l (t) - \delta_{kl} \overline{\chi}_l (t) \right), ~ \forall k.
\end{equation}

The goal of this paper is to develop an optimal bidding strategy for EV aggregators so that their profits are maximized in the electricity market. In addition, we hope that the optimal bidding strategy can also guarantee the maximization of social welfare, i.e., the overall utility of all agents in the market. The real-time bidding strategy as well as its properties is discussed in the following.

\section{Optimal Bidding Strategy of EV Aggregators}
\label{sec-III}
In this section, we focus on the development of the optimal bidding strategy. We first derive the strategy theoretically, and then a practical algorithm for implementation is proposed.

\subsection{Optimal Bidding Strategy}
To start, we compute the optimal solution of \textbf{P2}. For the sake of brevity, we denote
\bsq
\begin{align}
    \underline{w}_g(t) & \triangleq \frac{q_g(t) + z_g(t) - X_g}{V}, ~ g = 1, 2, \dots, G \\
    \overline{w}_g(t) & \triangleq \frac{q_g(t) + z_g(t)}{V}, ~ g = 1, 2, \dots, G.
\end{align}
\esq
Obviously, $\underline{w}_g(t) < \overline{w}_g(t)$. Then the optimal solution $\hat{x}_g(t)$ of \textbf{P2} is given by
\bsq
\begin{align} \label{eq:xg1}
    \hat{x}_g(t) = \left\{\begin{array}{ll}
      X_g(t), & \text{if}~ \pi(t) \in \left[0, \underline{w}_g(t) \right], \\
      V(\overline{w}_g(t) - \pi(t)), & \text{if}~ \pi(t) \in \left( \underline{w}_g(t), \overline{w}_g(t) \right],\\      
      0, & \text{if}~ \pi(t) \in \left( \overline{w}_g(t), +\infty \right),
    \end{array} \right.
\end{align}
or
\begin{align} \label{eq:xg2}
    \hat{x}_g(t) = \left\{\begin{array}{ll}
      V(\overline{w}_g(t) - \pi(t)), & \text{if}~ \pi(t) \in \left[0, \overline{w}_g(t) \right],\\      
      0, & \text{if}~ \pi(t) \in \left( \overline{w}_g(t), +\infty \right).
    \end{array} \right.
\end{align}
\esq
If $\underline{w}_g(t) > 0$, then $\hat{x}_g(t)$ is given by \eqref{eq:xg1}; otherwise it is given by \eqref{eq:xg2}. The results show that the optimal charging power of each group is a function of the electricity price in a specific time slot. Then, the total charging power of all groups in the time slot $t$ is given by
\begin{equation} \label{eq:x_def}
    \hat x(t) = \sum_{g = 1}^G \hat x_g(t), ~ \forall t.
\end{equation}
Suppose the EV aggregator chooses its charging power as specified in \eqref{eq:xg1} and \eqref{eq:xg2}. Then reversely, the electricity price of the aggregator can be also represented by its total charging power. Here, we define a function $h^{(t)}(\hat x)$ in the time slot $t$, which satisfies
\begin{equation} \label{eq:h_def}
    \pi (t) = h^{(t)} (\hat{x} (t)).
\end{equation}
Note that $h^{(t)}(\hat x)$ is the inverse function of the function given by \eqref{eq:xg1}, \eqref{eq:xg2}, and \eqref{eq:x_def}.
With these, the optimal bidding cost curve of the aggregator is set as a piecewise quadratic function below:
\begin{align} \label{eq:u_def}
    u^{(t)} (x) \triangleq \int_0^{x} h ^{(t)} (\xi) \text{d} \xi, ~ 0 \le x \le \sum_{g = 1}^G X_g.
\end{align}
The reason why the bidding cost curve takes the form of \eqref{eq:u_def} will be elaborated later by Proposition \ref{prp:agg}. Before progressing further, the concavity of $u^{(t)}(x)$ needs to be discussed as a prerequisite.
Denote the bidding cost curve in the time slot $t$ by $u^{(t)}(x)$, then we have the following proposition:
\begin{proposition} \label{prp:con}
    The optimal bidding cost curve $u^{(t)}(x)$ given by \eqref{eq:u_def} is concave on $\left[0, \sum_{g = 1}^G X_g(t) \right]$ for any $t$.
\end{proposition}
The proof of Proposition \ref{prp:con} can be found in Appendix \ref{appd:con}. The convexity of \eqref{obj:OPF} is then guaranteed by Proposition \ref{prp:con}, thus \textbf{P3} is a convex optimization problem.

Denote the charging power and electricity price of the aggregator $k$ in the optimal solution of \textbf{P3} in the time slot $t$ by ${x}_k^* (t)$ and $\pi_k^*(t)$, respectively, and denote the charging power of the group $g$ of the aggregator $k$ in the optimal solution of \textbf{P2} with given $\pi_k^*(t)$ by $\hat{x}_{kg} (t;\pi_k^*(t))$, then we have the following proposition:
\begin{proposition} \label{prp:agg}
If the optimal bidding cost curve is given by \eqref{eq:u_def}, and the electricity price in \eqref{obj:EV} is given by \eqref{eq:price}, then the EV charging schedule determined by social welfare maximization align with 
the interest of each price-taking EV aggregator, i.e.,
\begin{equation}
    {x}_k^* (t) = \sum_{g=1}^{G_k} \hat{x}_{kg} (t;\pi_k^*(t)), ~ g = 1, 2, \dots, G_k, ~ \forall k, ~ \forall t.
\end{equation}
where $G_k$ is the number of groups of the aggregator $k$.
\end{proposition}
The proof of Proposition \ref{prp:agg} can be found in Appendix \ref{appd:agg}. It is revealed by Propositions \ref{prp:con} and \ref{prp:agg} that the proposed bidding cost curve can maximize social welfare and the interests of aggregators simultaneously.

\subsection{Practical Algorithm to Generate the Bidding Cost Curve} \label{sec:online}
%In the time slot $t$, we solve a bilevel optimization problem, in which the upper-level and lower-level problems are \textbf{P2} and \textbf{P3} respectively. The optimal electricity prices $\hat{\pi}_k (t)$ of all aggregators obtained by solving \textbf{P3} are taken into \textbf{P2}, and the queues $q_g(t)$ and $z_g(t)$ of all aggregators are updated to derive the optimal bidding function in the time slot $t + 1$. The two procedures are repeated alternatively to complete the solution. In this section, the notation ``$(t)$" will be omitted for the sake of brevity since the algorithm can be applied to any time slot.

%\textbf{P2} can be easily solved as long as the electricity prices $\hat{\pi}_k(t)$ of all aggregators are given, while the solution of \textbf{P3} is not so straightforward because it is difficult to derive the analytical form of the optimal bidding function by \eqref{eq:u_def}. Instead, we calculate the integration by an online algorithm. In this section, we discuss how the online algorithm is developed and applied in this study.

%\subsection{Optimal Bidding Function}
In this section, we propose a practical algorithm (Algorithm 1) to generate the optimal bidding cost curve of an EV aggregator given in \eqref{eq:u_def}. Here, the notation ``$(t)$" will be omitted for the sake of brevity since the algorithm can be applied to any time slot. Moreover, since we focus on a single aggregator, we use the notation such as ``$u(x)$", in which the subscript ``$k$" is omitted for the sake of brevity. The strategies of other aggregators can be generated by the same approach.
From \eqref{eq:xg1} and \eqref{eq:xg2} we know that the optimal charging power of each group is a piecewise linear function. Thus, the optimal total charging power of all groups $\hat{x}$ is also a piecewise linear function.

We sort all elements of the sequence $(\underline{w}_g)_{g=1}^G$ and $(\overline{w}_g)_{g=1}^G$ in ascending order and denote the new sequence by $(\tilde{w}_g)_{g=1}^{2G}$, i.e., $\tilde{w}_1 \le \tilde{w}_2 \le \dots \le \tilde{w}_{2G}$. Here we assume $\tilde{w}_1 > 0$ without loss of generality, then $\hat{x}$ is supposed to take the following form:
%\ref{fig:bidding_curve}.
\begin{align} \label{eq:x_hat}
    \hat{x} = \left\{\begin{array}{ll}
      b_0, & \text{if}~ \pi \in [0, \tilde{w}_1],\\
      \dots & \\
      c_n \pi + b_n, & \text{if}~ \pi \in (\tilde{w}_{n-1}, \tilde{w}_n],\\
      \dots & \\
      0, & \text{if}~ \pi \in (\tilde{w}_{2G}, +\infty),
    \end{array} \right.
\end{align}
where $b_n$, $c_n$ are constants, $n = 0, 1, \dots, 2G$. Specifically, $b_0 = \sum_{g=1}^G X_g$. $c_n$ satisfies $c_n = -\beta_n V$, where $\beta_n$ is a non-negative integer. In the cases where there exists $\tilde{w}_n \le 0$, the segments defined over $[\tilde{w}_1, \tilde{w}_n]$ should be removed from \eqref{eq:x_hat}. 
$b_n$ and $c_n$ can be determined by a recursive method. Geographically, $c_n$ is the slope of the line segment in $(\tilde{w}_{n-1}, \tilde{w}_n]$. It can be observed from \eqref{eq:xg1} that the slope of $x_g$ falls by $V$ at $\underline{w}_g$ and rises by $V$ at $\overline{w}_g$ as $\pi$ is increased from 0. Therefore, $c_n, ~ n = 1, 2, \dots, 2G$ can be derived from an initial slope $c_0 = 0$ as follows:
\begin{equation}
    c_n = c_{n-1} + (\overline{\beta}_n - \underline{\beta}_n) V, ~ n = 1, 2, \dots, 2G,
\end{equation}
where $\underline{\beta}_n$, $\overline{\beta}_n$ are the numbers of elements equal to $c_n$ in $(\underline{w}_g)_{g=1}^G$ and $(\overline{w}_g)_{g=1}^G$, respectively. $b_n$ can be derived by the continuity of the piecewise linear function in \eqref{eq:x_hat}, i.e., the value of $\hat{x}$ at $\tilde{w}_{n-1}$ should both satisfy the expression of the section over $(\tilde{w}_{n-2}, \tilde{w}_{n-1}]$ and $(\tilde{w}_{n-1}, \tilde{w}_n]$, yielding
\begin{equation} \label{eq:b_n}
    b_n = b_{n-1} + (c_{n-1} - c_n) \tilde{w}_{n-1}, ~ n = 1, 2, \dots, 2G.
\end{equation}

Afterwards, the optimal bidding cost curve $u(x)$ can be obtained by piecewise integration. First we assume $\tilde{w}_1 \ge 0$, and denote the value of $\hat{x}$ at $\tilde{w}_n$ by $x_n, ~ n = 1, 2, \dots, 2G$. Specifically, let $x_0 = b_0$. Obviously, $c_n \neq 0$ when $x_n \neq x_{n-1}$. Since $\hat{x}$ is decreasing in its domain, $u(x)$ can be derived as follows:
\begin{equation} \label{eq:u_quad}
\begin{aligned}
    u(x) & = u(x_n) + \int_{x_n}^{x} \frac{\xi - b_n}{c_n} \text{d} \xi \\
     & = u(x_n) + \frac{x^2 - x_n^2}{2c_n} - \frac{b_n (x - x_n)}{c_n}, \\
     & x \in (x_n, x_{n-1}], \text{ if } x_n > x_{n-1}, ~ n = 1, 2, \dots, 2G.
\end{aligned}
\end{equation}
If there exists $\tilde{w}_n \le 0$, the domain of $u(x)$ shrinks to $[0, x'_0]$ because the electricity price $\pi$ must be non-negative, where $x'_0$ represents the value of $\hat{x}$ when $\pi = 0$.

The algorithm details are illustrated in Algorithm \ref{algo:bidding}, repeated in each time slot before solving \textbf{P3}. Note that the integration (steps \ref{state:int1} -- \ref{state:int2}) is performed backwards and stops when $\tilde{w}_{n-1} \le 0$.

\begin{algorithm}[ht]
\caption{Bidding Cost Curve Update} \label{algo:bidding}
\begin{algorithmic}[1]
\STATE Sort the elements of $(\underline{w}_g)_{g=1}^G$,  $(\overline{w}_g)_{g=1}^G$, and get $\{\tilde{w}_n\}$.
\STATE Initialize the slope and the intercept of the piecewise linear function: $c = 0, ~ b_0 = \sum_{g=1}^G X_g$.
\FOR{$n = 1:2G$}
\IF{$\tilde{w}_n$ is an element of $(\underline{w}_g)_{g=1}^G$}
\STATE $c = c - V$
\ELSE
\STATE $c = c + V$
\ENDIF
\STATE $c_n = c$, compute $b_n$ by \eqref{eq:b_n}.
\ENDFOR
\FOR{$n = 2G:-1:1$} \label{state:int1}
\IF{$c_n \neq 0$}
\STATE Compute $u(x)$ in $(x_n, x_{n-1}]$ by \eqref{eq:u_quad}.
\ENDIF
\IF{$\tilde{w}_{n-1} \le 0$}
\STATE Break the loop.
\ENDIF
\ENDFOR \label{state:int2}
\end{algorithmic}
\end{algorithm}

\subsection{Overall Market Operation}
\begin{figure}
    \centering
    \includegraphics[width=\linewidth]{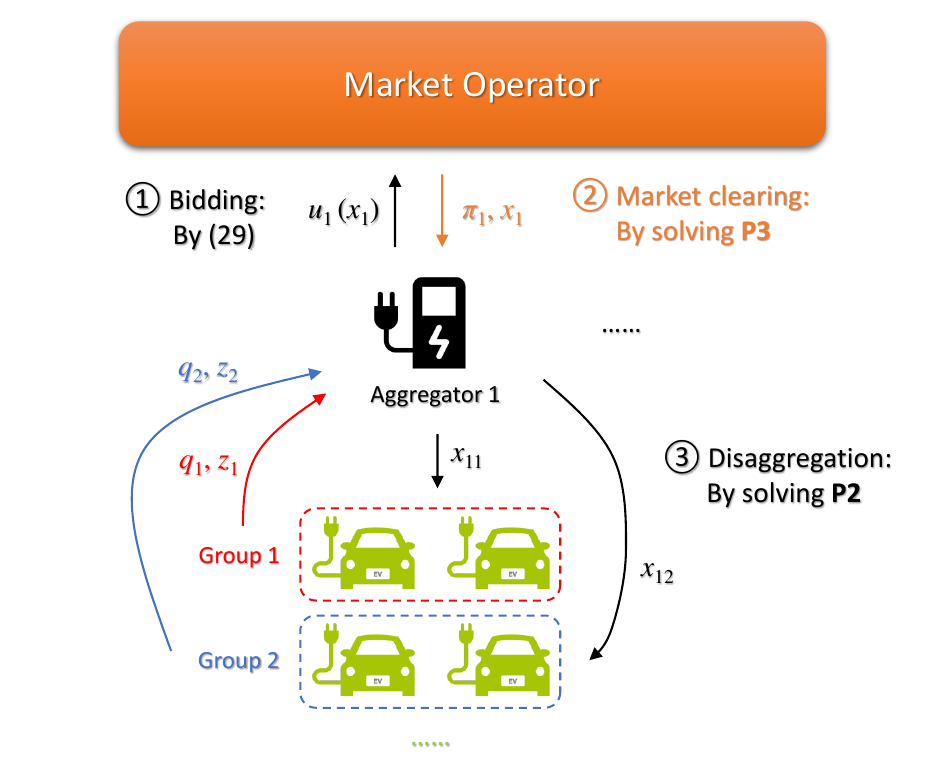}
    \caption{The optimal bidding framework.}
    \label{fig:framework}
\end{figure}

So far, we have established the real-time optimal bidding framework for EV aggregators. In each time slot $t$, an aggregator bids the function $u^{(t)} (x)$ to the market operator based on its backlogs of queues. Upon receiving the bidding cost curves from all aggregators, the operator solves \textbf{P3} to clear the market, by which it determines the electricity prices of aggregators. Then the operator returns the price to each aggregator. Afterwards, each aggregator solves \textbf{P2} using the latest price, and update the charging power of EVs in each group. The procedure is repeated in each time slot. An illustration of the overall "aggregator bidding-market clearing" framework is presented in Fig. \ref{fig:framework}.

However, with the aggregators' bidding cost curves by \eqref{eq:u_quad}, the market clearing problem \textbf{P3} cannot be solved directly by commercial software since the bidding cost curves are piecewise quadratic functions. In the following, we introduce a linearization technique, by which the approximate solution of \textbf{P3} can be obtained. Thereafter, the objective \eqref{obj:OPF} can be simplified and thus the accurate solution of \textbf{P3} can be finally procured.

\textbf{P3} can be approximated by the following problem:
\bsq
\begin{align}
    \textbf{P4:}~ \min_{p_i, x_k, y_k}~ & \sum_{i = 1}^I f_i^{\text{Gen}}(p_i) - \sum_{k = 1}^K y_k \\
    \mbox{s.t.}~ & x_k = \sum_{m=1}^{M_k} \sigma_{km} x_{km}, ~ \forall k, \\
    & y_k = \sum_{m=1}^{M_k} \sigma_{km} u_k (x_{km}), ~ \forall k, \\
    & \sigma_{km} \ge 0, ~ m = 1, 2, \dots, M_k, ~ \forall k, \\
    & \sum_{m=1}^{M_k} \sigma_{km} = 1, ~ \forall k, \eqref{eq:ErgBal} - \eqref{ineq:LinePwr},
\end{align}
\esq
where $M_k$ is a large integer, $\forall k$. $(x_{km}, u_k (x_{km})), ~ m = 1, 2, \dots, M_k$ are the sample points of $u_k(x)$, and $(x_k, y_k)$ is the convex combination of these points. 
The optimal solution of \textbf{P4} will be sufficiently close to the accurate solution of \textbf{P3} as long as $M_k$ is large enough. However, there is a trade-off as the computational cost will also be increased if a larger $M_k$ is chosen. Under the circumstances where high accuracy is not required, the approximate optimal solution should be good enough. Alternatively, we derive the accurate solution of \textbf{P3} using the solution of \textbf{P4} in this study. Derived by \eqref{eq:u_quad}, the optimal bidding cost curve of any aggregator $k$ $u_k(x_k)$ is a piecewise quadratic function. 
The interval in which $\hat{x}_k$ falls can be revealed by $\bar{x}_k$ if $M_k$ is large enough, thus $u_k(x_k)$ in \eqref{obj:OPF} can be substituted by a quadratic function.
Suppose $\hat{x}_k \in (x_{k n_k}, x_{k n_k - 1}]$, then \textbf{P3} is equivalent to
\bsq
\begin{align}
    \textbf{P5: } \min_{p_i, x_k}~ & \sum_{i = 1}^I f_i^{\text{Gen}}(p_i) - \sum_{k = 1}^K u_{k n_k}(x_k) \label{obj:accr} \\
    \mbox{s.t.}~ & x_k \in (x_{k n_k}, x_{k n_k - 1}], ~ \forall k, \\
    & \eqref{eq:ErgBal}, ~ \eqref{ineq:GenPwr}, ~ \eqref{ineq:LinePwr}, \nonumber
\end{align}
\esq
where $u_{k n_k}(x)$ is the optimal bidding cost curve of the aggregator $k$ in $(x_{k n_k}, x_{k n_k - 1}]$.
\textbf{P5} is a quadratic problem with linear constraints and thus can be solved by its KKT condition. The details of the "aggregator bidding-market clearing" framework is presented in Algorithm \ref{algo:quadratic}. %\jpang{So far, we have been proving everything for quadratic bidding cost curves. Perhaps want to state something around extension? }.

\begin{algorithm}[ht]
\caption{Aggregator Bidding-Market Clearing} \label{algo:quadratic}
\begin{algorithmic}[1] 
\STATE Set up the case and initialize the variables $x_k, ~ \forall k$.
\FOR{$t = 1 : T$}
\FOR{$k = 1 : K$}
\STATE Run Algorithm \ref{algo:bidding} and collect $u_{k1}(x)$, $u_{k2}(x)$, $\dots$, $u_{k n_k}(x)$.
\ENDFOR
\STATE Solve \textbf{P4}.
\FOR{$k = 1 : K$}
\FOR{$n = 1 : 2G_k$}
\IF{$\hat{x}_k \in (x_{k n_k}, x_{k n_k - 1}]$}
\STATE Take $u_{k n_k}(x)$ into \eqref{obj:accr} and break the loop.
\ENDIF
\ENDFOR
\ENDFOR
\STATE Solve \textbf{P5}.
\STATE Update the electricity price of each aggregator by \eqref{eq:price}.
\STATE Solve \textbf{P2}.
\STATE Update $q_g$ and $z_g$ by \eqref{eq:Qgt} and \eqref{eq:zgt} respectively.
\ENDFOR
\end{algorithmic}
\end{algorithm}

\section{Case Studies} \label{sec:result}

In this section, the proposed algorithm is first validated and its performance is evaluated against other approaches afterward. 
We also examine the influence of several parameters on the results.
Finally, the scalability of the proposed algorithm is examined.

\subsection{Simulation Setup} \label{subsec:setup}
The investigated system is programmed in MATLAB 2022a and the proposed algorithm is implemented on a standard desktop PC with an Intel i5-10505 CPU and 16 GB RAM.
Case studies in this paper are performed in an IEEE 118-bus network with 54 generators and 23 fixed loads. The investigated time horizon is 24 hours, with 5 minutes for each time slot. 
We model a total of 10 EV aggregators in the network, covering 3 types of maximum charging power: 7 kW, 60 kW and 120 kW, depending on the charging facility of each aggregator, and the charging efficiency $\eta = 0.95$ in all cases.
EVs in the service area are assumed to have one of 3 types of battery capacity: 20 kWh, 50 kWh and 80 kWh, which are mainstream specifications in the EV market. The state of charge (SOC), i.e., the percentage of battery energy, of each EV upon departure is 80\% and the SOC on arrival is set between 30\%--70\%. 
The EVs served by each aggregator are divided into 10 groups according to their duration of parking. The size of each group, i.e., the number of EVs in the group, ranges between 30--50. 
% Moreover, EVs are assumed to arrived at the aggregators at a random time in the morning and leave in the evening, as shown in Fig. \ref{fig:arrival}. Most EVs arrive between 7--10 am and leave between 18--22 pm. Hence, the above settings are able to simulate the real-world behaviors of EVs in an urban or industrial zone to some extent.

% \begin{figure}[ht]
%   \centering
%   \includegraphics[width = 0.4\textwidth]{Figs/arrivals.eps} \\
%   \caption{Histogram of time of arrival and departure at an EV aggregator.} \label{fig:arrival}
% \end{figure}

\subsection{Validation of the Proposed Algorithm}
\begin{figure}[htbp]
    \centering
    \subfigure[]
    {\includegraphics[width=4cm]{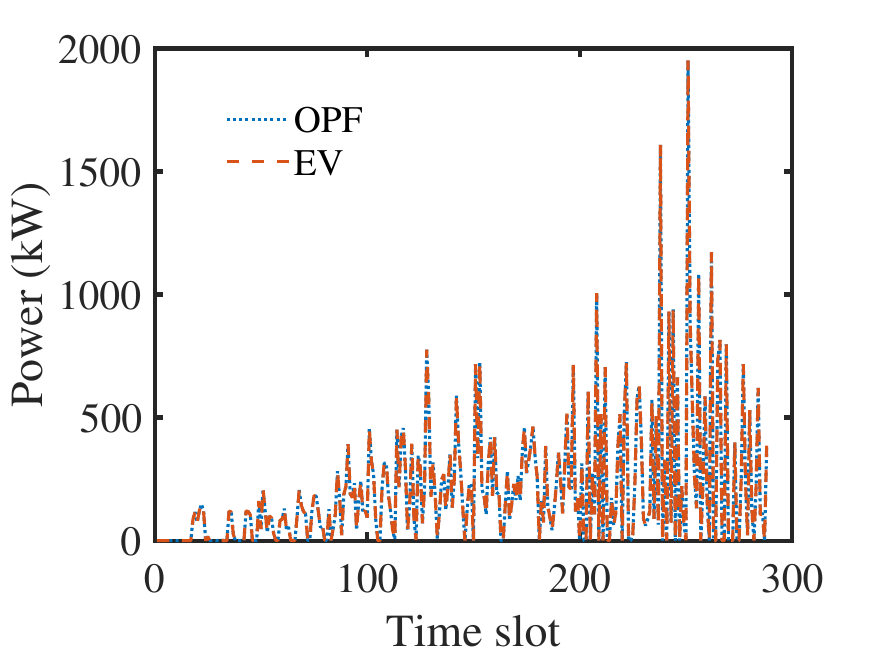} \label{fig:vld}}
    \subfigure[]
    {\includegraphics[width=4cm]{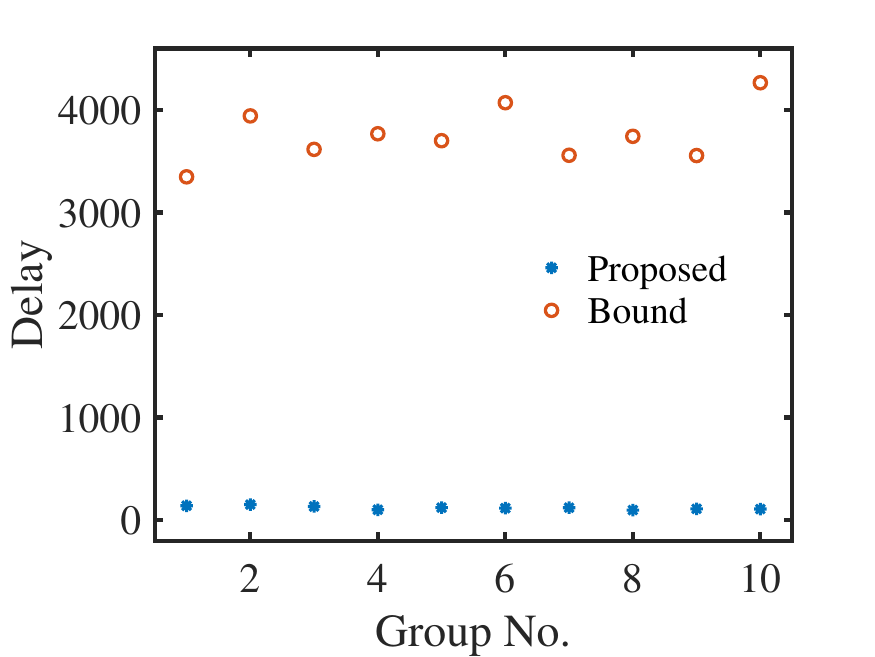} \label{fig:delay_g}}
    \caption{(a) The optimal charging power of an aggregator by solving \textbf{P2} (EV) and \textbf{P3} (OPF). (b) The maximum charging delay of all groups of an aggregator and the theoretical bounds.}
\end{figure}

First of all, we verify the correctness of the proposed algorithm. According to Proposition \ref{prp:agg}, the charging power of each aggregator in the optimal solution of \textbf{P2} should be consistent with the total of all groups of the aggregator in the optimal solution of \textbf{P3}. Part of the results are shown in Fig. \ref{fig:vld}, in which EV and OPF refer to the optimal solutions of \textbf{P2} and \textbf{P3} respectively. It can be observed that the power to be disaggregated by the aggregator, which is assigned by the central operator of the network, perfectly matches the total charging power of all groups of the aggregator. Hence, Proposition \ref{prp:agg} is validated.

The maximum charging delay of all groups of EVs served by an aggregator is presented in Fig. \ref{fig:delay_g} as well as the theoretical bounds given by \eqref{ineq:delay}. The results show that the maximum charging delay obtained by the proposed algorithm is far below the bound for each group. Through this, we validate Proposition \ref{prp:delay}.

% \begin{figure}[ht]
%   \centering
%   \includegraphics[width = 0.45\textwidth]{Figs/validation.eps} \\
%   \caption{The optimal charging power of 4 aggregators in the network obtained by solving \textbf{P2} and \textbf{P3}.} \label{fig:vld}
% \end{figure}

\begin{figure}[htbp]
    \centering
    \subfigure[]
    {\includegraphics[width=4cm]{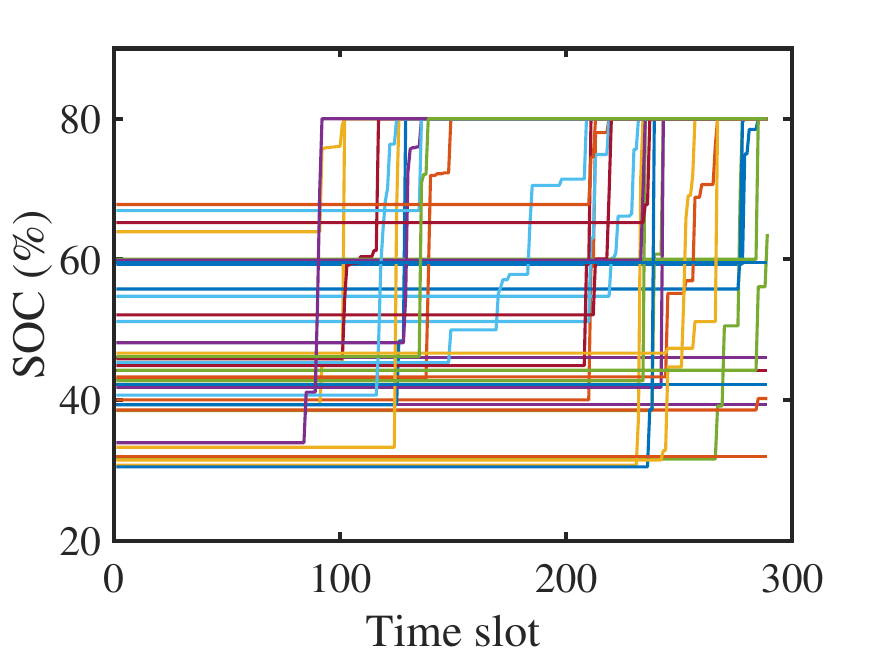} \label{fig:soc}}
    % \subfigure[]
    % {\includegraphics[width=4cm]{Figs/roa_ag1.eps} \label{fig:roa}}
    \subfigure[]
    {\includegraphics[width = 4cm]{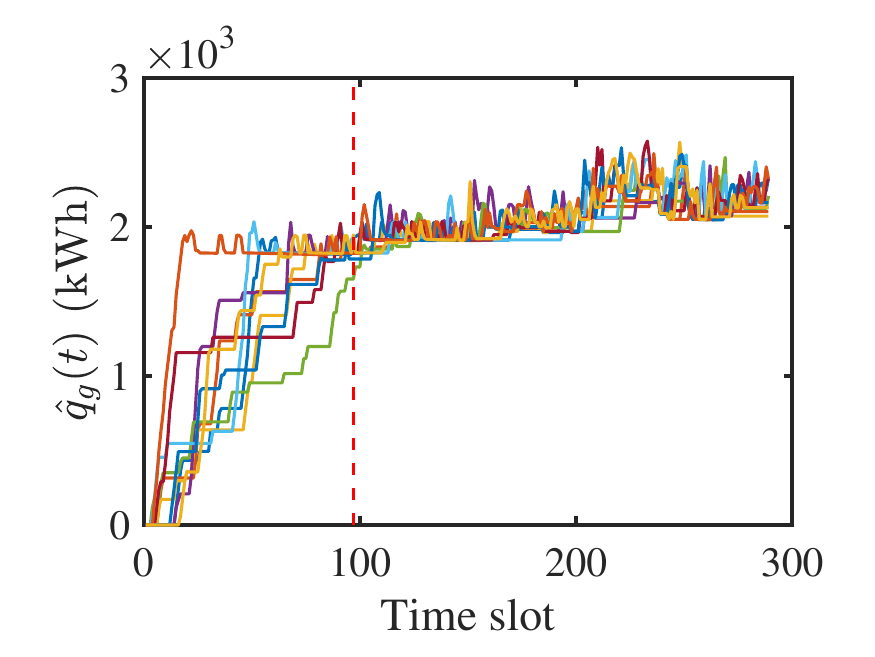} \label{fig:q}}
    \subfigure[]
    {\includegraphics[width = 4cm]{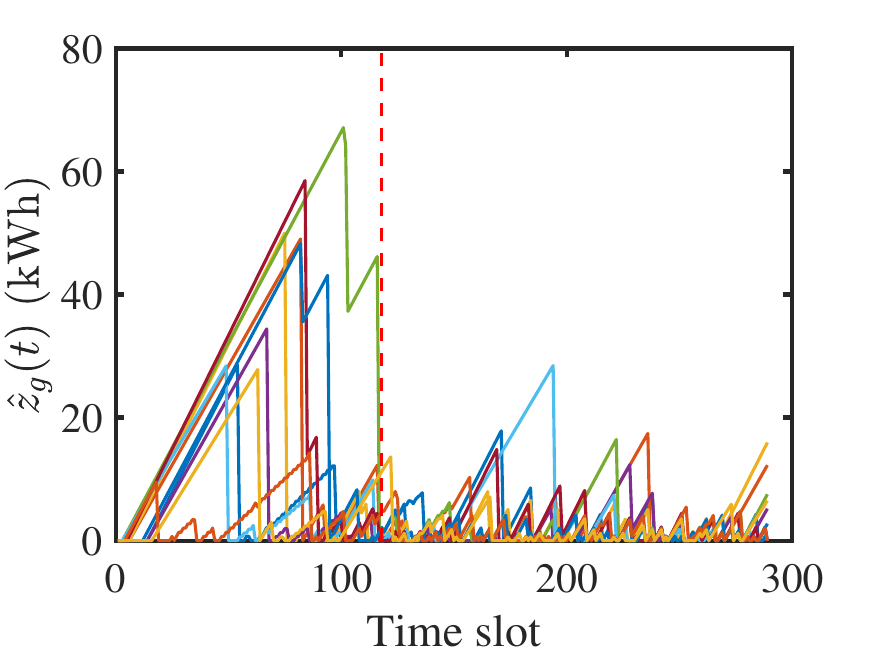} \label{fig:z}}
    % \subfigure[]
    % {\includegraphics[width = 4cm]{Figs/xg_hat_ag1.eps} \label{fig:xg}}
    \subfigure[]
    {\includegraphics[width = 4cm]{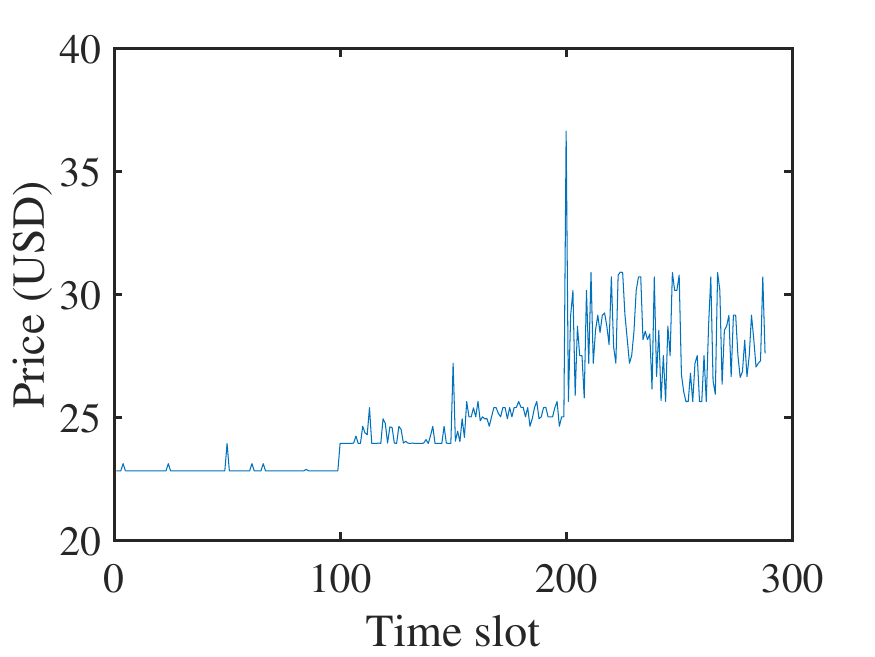} \label{fig:price}}
    \caption{Results of simulation for one aggregator: (a) The SOC traces of a group of EVs; (b) The queues of charging demands; (c) The virtual queues; (d) The electricity price.} \label{fig:agg}
\end{figure}

Subsequently, we need to check whether all the charging demands have been met at the end of the time horizon. The results related to the quality of service of an aggregator are presented in Fig. \ref{fig:agg}. 
The progress of SOC of the EVs in a group from arrival to departure is shown in Fig. \ref{fig:soc}. There are 40 EVs in the group, and the SOC of each EV on arrival ranges from 30\%--68\%.
% For each group of the aggregator, the EVs are charged at high power on their arrival and most of them are fully charged within xxx time slots except one, which is charged at lower power in the remaining time. 
By the end of the time horizon, most of the EVs have been fully charged, i.e., SOC $\ge$ 80\%. There are still 8 EVs below the desired energy level since they arrive at the charging station later than the others. As a result, even using an offline model, these EVs cannot be fully charged at the end of the time horizon studied.
% {\color{blue}However, based on the performance in the investigated period, it is reasonable to assume that they can be also fully charged by the time they leave the charging station.} 
% The percentage of EVs fully charged for each group is shown in Fig. \ref{fig:roa}, and mostly ranging from 60\%--80\%, which is a medium level of quality of service. 
In this case, we set $V = 80$. The relationship between $V$ and the percentage of EVs fully charged will be elaborated in Section \ref{subsec:V}.

The queues of charging demands $q_g(t)$ for all groups of the aggregator are presented in Fig. \ref{fig:q}. 
% corresponding to the high-power charging and low-power charging phases observed in Fig. \ref{fig:soc}. 
At the beginning of the time horizon, there is no delayed charging demand, and the backlogs of $q_g(t)$ are 0 for all groups. Two phases can be recognized in the behavior of $q_g(t)$ for each group. In the first phase, $q_g(t)$ increases sharply from 0. Since the backlogs are rather low and charging delay is allowed, the cost minimization dominates over the Lyapunov drift in the objective function of \textbf{P2}. As the backlogs of $q_g(t)$ increase, the stability of queues becomes more important in the second phase, thus the charging power of each group is correspondingly increased to consume the charging demands. $q_g(t)$ remains nearly constant at the end of the time horizon.
% It can be observed that most of the charging demands are processed in a short period after arrival and then the queues are consumed at a lower rate in the second phase since charging delay is allowed. Finally, $q_g(t)$ diminishes to 0 for each group, which indicates that all the charging demands have been satisfied. 
Similarly, the development of the virtual queues $z_g(t)$ of all aggregators can be also divided into two phases, as shown in Fig. \ref{fig:z}. In the first phase, the backlogs of $z_g(t)$ for groups rise sharply as EVs arrive, while all of them fluctuate around 0 in the second phase due to the large availability of charging power.
% It finally remains constant after $q_g(t)$ is decreased to 0, indicating the accomplishment of charging service for all EVs.
The electricity price that the aggregator accepts during the investigated period is shown in Fig. \ref{fig:price}. In our simulation, the electricity price is influenced by both the charging demands and other loads in the power grid. While the price fluctuates smoothly between 23--27 USD/kWh before the time slot 199, the time slot 200 has seen a dramatic surge and a turbulent variance arises afterward due to the increase of other loads in our settings.

\subsection{Performance of the Proposed Algorithm}

We further compare the results of the proposed algorithm against the following benchmarks to evaluate its performance:
\begin{enumerate} [ B1.]
    \item The offline method, i.e., solving \textbf{P1}; \label{bm:offline}
    \item A linear counterpart of the proposed algorithm introduced in Appendix \ref{appd:lin}; \label{bm:lin}
    \item A simplified online algorithm with alternative bidding. \label{bm:simp}
\end{enumerate}
The idea of B\ref{bm:simp} is to replace the utility function $u_{kt}(x_k(t))$ in the bidding by the bounds of $x_k(t)$, that is, \textbf{P3} is simplified as the following problem:
\bsq
\begin{align}
    \textbf{P3$''$: } \min_{p_i, x_k}~ & \sum_{i = 1}^I f_i^{\text{Gen}}(p_i(t)) \\
    \mbox{s.t.}~ & \underline{X}_k(t) \leq x_k(t) \leq X_k(t), ~ \forall k, ~ \forall t, \label{ineq:x_bd_new} \\
    & \eqref{eq:ErgBal}, ~ \eqref{ineq:GenPwr}, ~ \eqref{ineq:LinePwr}, \nonumber
\end{align}
\esq
where a new lower bound of $x_k(t)$, denoted by $\underline{X}_k(t)$, is introduced into the constraint \eqref{ineq:x_bd_new}. $\underline{X}_k$ is defined as follows:
\bsq
\begin{equation} \label{eq:Xlb}
    \underline{X}_k(t) = \max \left\{\underline{X}_k^{(1)}(t),~ \underline{X}_k^{(2)}(t) \right\},
\end{equation}
where
\begin{align}
    \underline{X}_k^{(1)}(t) & = \sum_{g=1}^{G_k} \sum_{v \in \mathcal{S}_g} \frac{E_v^{\text{d}} - e_v(t)}{\eta (T_v^{\text{d}} - t)}, \label{eq:Xlb1} \\
    \underline{X}_k^{(2)}(t) & = \sum_{g=1}^{G_k} \sum_{v \in \mathcal{S}_g} \frac{E_v^{\text{d}} - e_v(t) - \eta P_v (T_v^{\text{d}} - t - 1)}{\eta \Delta t}. \label{eq:Xlb2}
\end{align}
\esq
\eqref{eq:Xlb1} and \eqref{eq:Xlb2} are actually two ways of defining the lower bound of $x_k(t)$. While the former chooses the average power calculated by the rest of energy and time, the latter expects that the charging task should be finished immediately assuming the EV has been charged at the maximum power since the beginning of charging. We adopt a hybrid definition in \eqref{eq:Xlb} so that $\underline{X}_k(t)$ will not be too low in the middle of charging or too high in the late stage.

\begin{figure}[htbp]
    \centering
    \includegraphics[width = 0.6\linewidth]{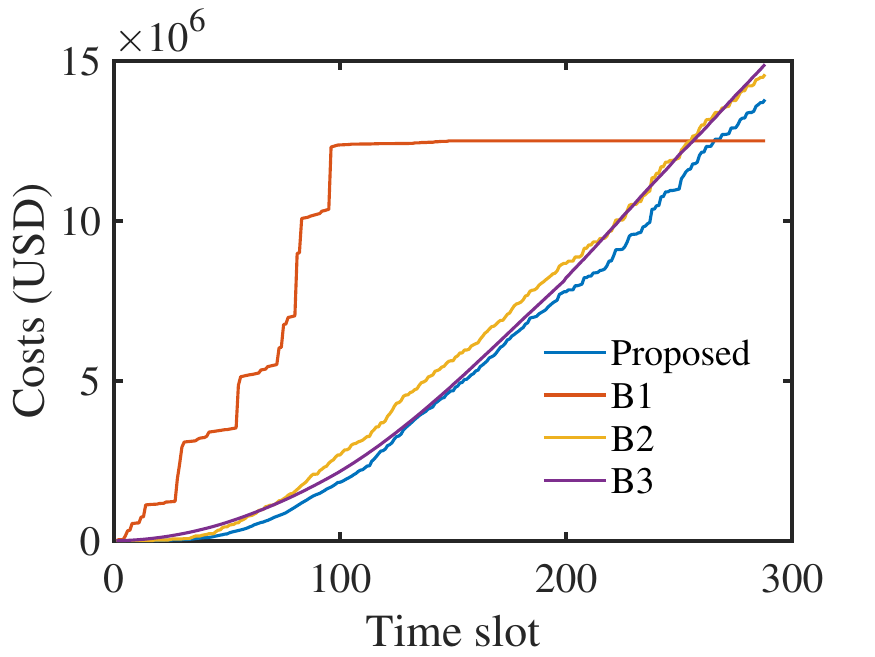} \label{fig:costs}
    \caption{A comparison between the proposed algorithm and the benchmarks.}
\end{figure}

To evaluate the performance of the proposed algorithm, we compare the accumulated costs, i.e., the sum of the costs from the beginning of the time horizon to the time slot $t$, obtained by the proposed algorithm, and the benchmarks listed above. The traces of the accumulated costs against time are shown in Fig. \ref{fig:costs}. Obviously, the offline benchmark B\ref{bm:offline} has the lowest total costs, although it is not always the best during the investigated period. Nevertheless, the performance of the proposed algorithm does not fall far from B\ref{bm:offline}. The total costs of all algorithms are listed in Table \ref{tab:costCmp}, and the relative values are also given with B\ref{bm:offline} as the benchmark.
% The results show that the total costs of the proposed algorithm is 5\% higher than B\ref{bm:offline}. 
The unit cost is the ratio of the total cost to the total purchase of electricity.
%, while the optimality gap given by Proposition \ref{prp:gap} is {\color{red} xxx\%}.
In this case, the left-hand and the right-hand side of \eqref{ineq:gap} are $4.48 \times 10^3$ USD and $7.59 \times 10^5$ USD respectively. Hence, Proposition \ref{prp:gap} is validated.
It can be also observed that the proposed algorithm has a remarkable advantage over its linear counterpart B\ref{bm:lin}, with a reduction of 5\% in the total costs. This advantage stems from the quadratic term in the objective of the proposed algorithm \eqref{obj:EV}, which is relaxed to its upper bound in the derivation of the linear counterpart. The benchmark B\ref{bm:simp} has the lowest costs in most of the timesteps (even better than B\ref{bm:offline}). However, the costs of B\ref{bm:simp} are continuously rising at the end of the time horizon, leading to the highest final outcome among all the algorithms. In summary, the proposed algorithm has a competitive performance compared with the benchmarks, and the gap between the proposed algorithm and the offline optimum is acceptable.
% Fig. \ref{fig:roa} shows the completion rate of the charging task, i.e., the percentage of EVs that have been charged to 80\%, in each group for an aggregator. It can be observed that the proposed algorithm has the highest completion rate among all algorithms, more than twice of that of B\ref{bm:offline}. Therefore, the proposed algorithm is able to maintain a high level of quality of service with limited cost.

% \begin{figure}[htbp]
%     \centering
%     \subfigure[]
%     {\includegraphics[width = 4cm]{Figs/costs.eps} \label{fig:costs}}
%     \subfigure[]
%     {\includegraphics[width = 4cm]{Figs/roa_ag7.eps} \label{fig:roa}}
%     \caption{(a) A comparison between the proposed algorithm and the benchmarks. (b) Completion rate of charging task.}
% \end{figure}

\begin{table}[!htbp]
    \centering
    \caption{Cost comparison between algorithms.} \label{tab:costCmp}
    \begin{tabular}{@{}ccccccccc@{}}
        \toprule
        \multirow{2}{*}{Algorithm} & Total cost & \multirow{2}{*}{Relative value} & Unit cost \\
        &  (Million USD) & & (USD/kWh) \\
        \midrule
        Proposed & 13.79 & 110\% & 25.21 \\
        B1 & 12.50 & 100\% & 22.85 \\
        B2 & 14.58 & 117\% & 25.08 \\
        B3 & 14.89 & 119\% & 25.95 \\        
        \bottomrule
   \end{tabular}
\end{table}

\subsection{The Impact of $V$ and $\alpha_g$} \label{subsec:V}

\begin{figure}[htbp]
    \centering
    \subfigure[]
    {\includegraphics[width=4cm]{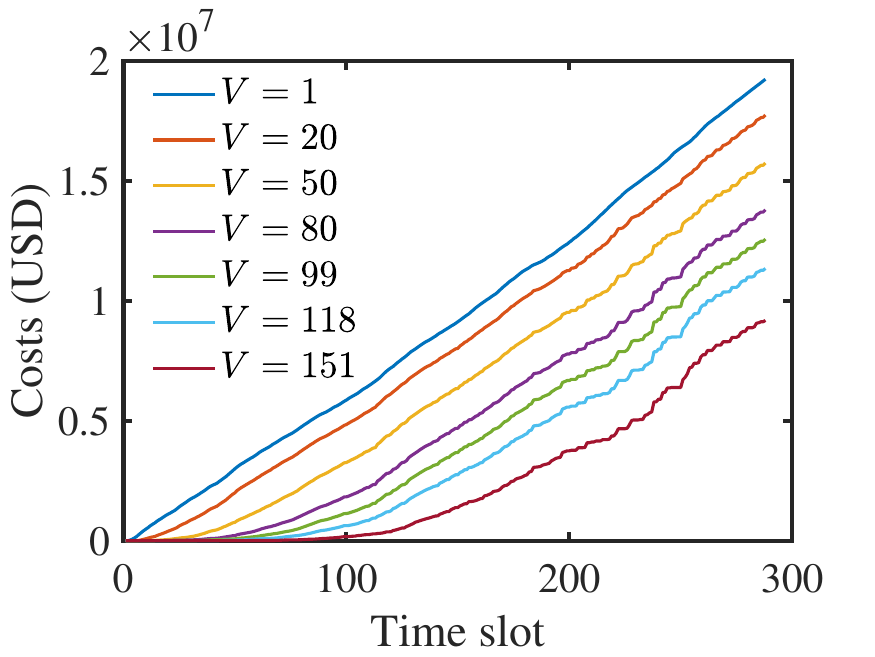} \label{fig:costs_v}}
    % \subfigure[]
    % {\includegraphics[width=4cm]{Figs/roa_v_80_ag7.eps} \label{fig:roa_v_80_ag7}}
    \subfigure[]
    {\includegraphics[width=4cm]{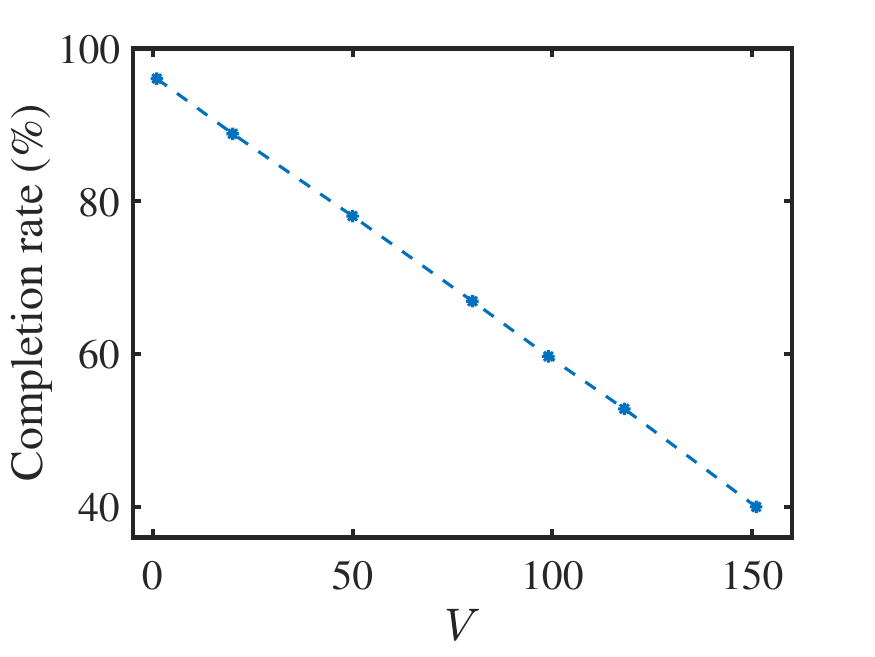} \label{fig:roat_v}}
    \subfigure[]
    {\includegraphics[width=4cm]{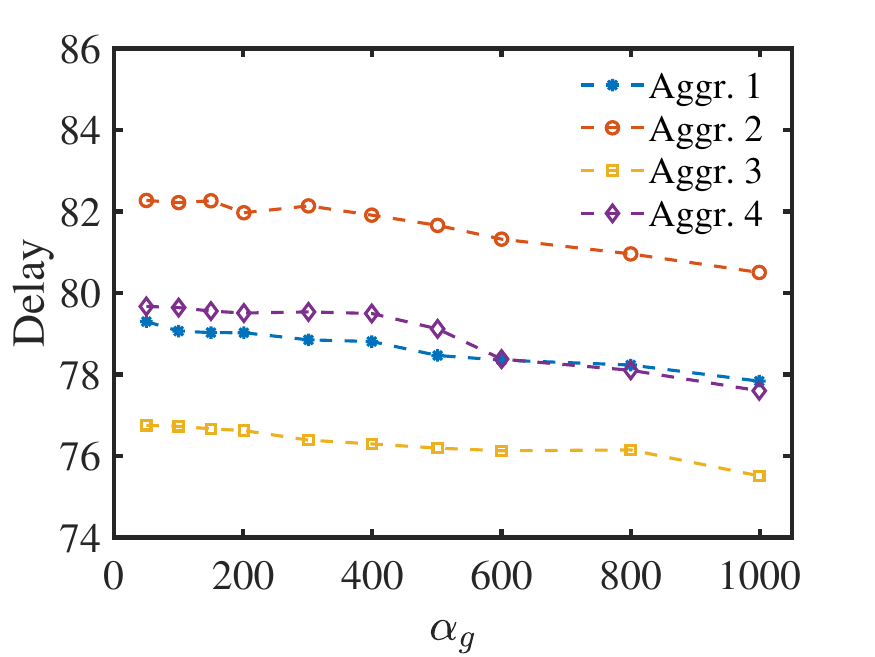} \label{fig:delay_ag}}
    \subfigure[]
    {\includegraphics[width=4cm]{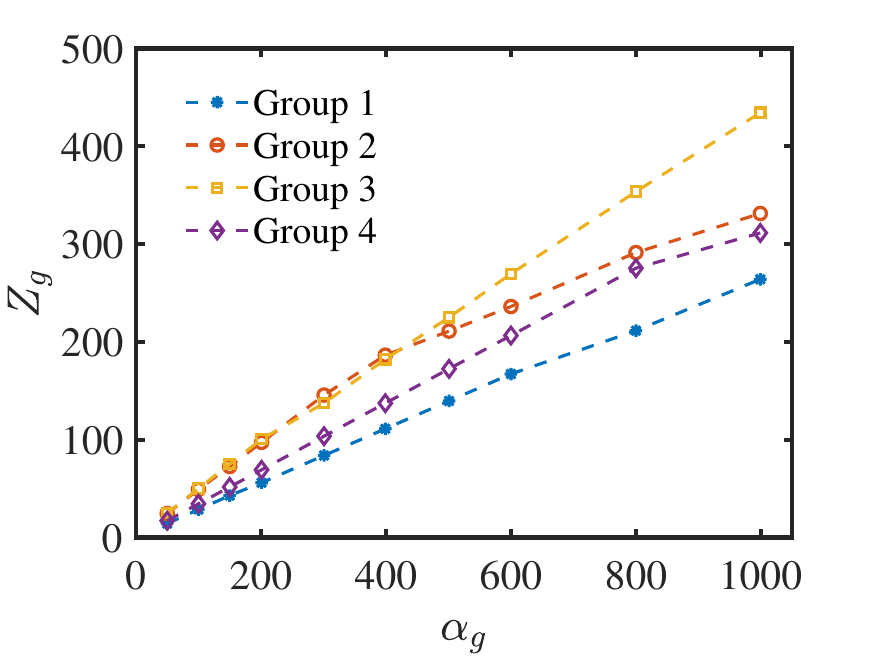} \label{fig:zg_alpha}}
    \caption{(a) The total costs of all aggregators with different $V$. (b) The overall completion rate of all aggregators. (c) The average charging delay of EVs for some aggregators. (d) The maximum backlogs of virtual queues for some groups of an aggregator.}
\end{figure}

As mentioned in Section \ref{subsec:lya}, the constant $V$ is a parameter balancing the stability of queues and optimality. In this section, the impact of $V$ on the performance of the proposed algorithm is investigated. The accumulated costs of all EV aggregators in the model against time are compared with $V$ ranging from 1 to 151. 
As shown in Fig. \ref{fig:costs_v}, the cost decreases when $V$ increases, caused by the gaining weight of the costs $f(t)$ in the objective function \eqref{obj:driftPlusP} against the stability of queues.
However, the reduction of operational costs of aggregators comes at a price of deteriorating quality of service. It can be observed in Fig. \ref{fig:roat_v} that the completion rate of charging tasks for EVs served by all aggregators decreases as $V$ increases. Specifically, there is an approximately linear relationship between the completion rate and $V$.
%increasing delay of EV charging. 

As indicated in the definition of $z_g (t)$, $\alpha_g$ is a constant that controls the behavior of virtual queues. It can be deduced from \eqref{eq:zgt} that the backlogs of $z_g (t)$ tend to add up more quickly with larger $\alpha_g$, thus the remaining charging demands will be satisfied more quickly. To evaluate the impact of $\alpha_g$ on charging delay, the average charging delays of 4 aggregators are presented in Fig. \ref{fig:delay_ag} for clarity. It shows that the charging delay indeed decreases as $\alpha_g$ increases despite slight fluctuation, which is caused by the increasing $Z_g$, i.e., the maximum backlogs of virtual queues, shown in Fig. \ref{fig:zg_alpha}.
% In an extreme case, the charging demands cannot be met when $V$ exceeds a certain value. 
% Denote this threshold by $V_{\text{cr}}$. In this study, $V_{\text{cr}}$ falls between 1 and 3 for each aggregator. 
% It can be observed from Fig. \ref{fig:costs_v} that the cost of aggregators when $V = 3$ is remarkably lower than those with smaller $V$, caused by unfinished charging task. For example, the SOC traces of an EV group of an aggregator is shown in Fig. \ref{fig:soc_v3}. Most of the EVs in the group have fulfilled the minimum energy level (SOC $\ge$ 80\%) at departure, whereas one of them fails to complete the charging task.

% \begin{figure}[ht]
%   \centering
%   \includegraphics[width = 0.45\textwidth]{Figs/costs_v.eps} \\
%   \caption{The total costs of all aggregators with different $V$.} \label{fig:costs_v}
% \end{figure}

% \begin{figure}[ht]
%   \centering
%   \includegraphics[width = 0.45\textwidth]{Figs/soc_v3_ag8_g1.eps} \\
%   \caption{The progress of SOC of the EVs in a group, $V = 3$.} \label{fig:soc_v3}
% \end{figure}

\subsection{Scalability}

The size of real-world EV optimal bidding networks can be much larger than the numerical experiments above. As such, the performance of proposed algorithm against a growing problem size, especially its computational burdens, also needs to be investigated. To evaluate the scalability of proposed algorithm, two groups of tests are performed by expanding the model size: 1) Increasing the number of buses in the network while the number of aggregators remains; 2) Increasing the number of aggregators while the number of buses remains. Each case is tested 50 times to reduce random errors. The results are shown in Table \ref{tab:scal}. Note that the specific computational time, i.e., the ratio of the computational time to the number of buses or aggregators, is presented here instead of the computational time, because the former is more meaningful when it comes to the efficiency of an algorithm on models of different sizes. The results show that the specific computational time with respect to whether the number of buses or that of aggregators is decreasing, which means the growing speed of computational time is slower than that of the scale of the problem. Hence, one could observe that the scalable method could easily adapt to the demands of large networks.

% \begin{figure}[htbp]
%     \centering
%     \subfigure[]
%     {\includegraphics[width=4cm]{Figs/scal.eps} \label{fig:scal}}
%     \subfigure[]
%     {\includegraphics[width=4cm]{Figs/scal_ag.eps} \label{fig:scal_ag}}
%     \caption{The results of scalability tests: (a) Increasing number of buses; (b) Increasing number of aggregators.}
% \end{figure}

\begin{table}
    \centering
\caption{Results of scalability tests}
\label{tab:scal}
    \begin{tabular}{cccc}
    \toprule
         Number of buses&  30&  57& 118\\ 
         Specific computational time (s)& 2.22 & 1.12 & 0.61 \\ 
         \midrule
         Number of aggregators&  5&  10& 20\\ 
         Specific computational time (s)& 11.14 & 7.24 & 4.99 \\ 
         \bottomrule
    \end{tabular}  
\end{table}

\section{Conclusion} \label{sec:conclu}
In this paper, the optimal real-time bidding strategy for EV aggregators is developed based on Lyapunov optimization techniques. Several nice properties of the proposed bidding strategy are proven theoretically and numerically. In addition, a practical algorithm based on stepwise integration is proposed to generate the bidding cost curve in a programmatic manner. Numerical experiments validate the proposed bidding strategy and reveal the following findings:
\begin{enumerate}
    \item The proposed bidding strategy can maximize the profits of EV aggregators and the social welfare simultaneously.
    \item The proposed real-time bidding strategy outperforms other online bidding strategies in terms of total cost and completion rate of charging task.
\end{enumerate}

Future efforts will be devoted to considering more complicated scenarios, such as the difference between aggregators on the convenience of service and EV preferences.

% The results of case studies show that the proposed method is practical and efficient comparing with selected benchmarks. The charging demands of all EVs in the simulation have been satisfied before their departure. To evaluate the performance of the proposed method, we compare the cost of buying electricity of all EV aggregators. The gap between the proposed method and the offline benchmark is less than the theoretical bound, and the cost of the proposed method is lower than the online benchmark. 

\bibliographystyle{IEEEtran}
\bibliography{ref}

\clearpage

\appendices
\makeatletter
\@addtoreset{equation}{section}
\@addtoreset{theorem}{section}
\makeatother
\setcounter{equation}{0}  
\renewcommand{\theequation}{A.\arabic{equation}}
\renewcommand{\thetheorem}{A.\arabic{theorem}}

\section{Proof of Proposition \ref{prp:gap}} \label{appd:gap}

\begin{proof}
Taking the optimal solutions of \textbf{P1} and \textbf{P2} into \eqref{ineq:dpp}, we have
\begin{equation} \label{ineq:cmp_f1}
    \begin{aligned}
        & \mathbb{E} \left[ \Delta(\boldsymbol{\Theta} (t)) + V \hat{f} (t) \right] \\
        \le ~& M + V \mathbb{E} \left[ \hat{f} (t) \right] - \sum_{g = 1}^G \mathbb{E} \left[ (q_g(t) + z_g(t)) \hat{x}_g(t) \right] \\
         & + \frac{1}{2} \sum_{g = 1}^G \mathbb{E} \left[ \hat{x}_g^2(t) \right] \\
         \le ~& M + V \mathbb{E} \left[f^*(t)\right] - \sum_{g = 1}^G \mathbb{E} \left[ (q_g(t) + z_g(t)) x^*_g(t) \right] \\
         & + \frac{1}{2} \sum_{g = 1}^G \mathbb{E} \left[ (x_g^*(t))^2 \right], ~ t = 1, 2, \dots, T.
    \end{aligned}
\end{equation}
Replace $\Delta(\boldsymbol{\Theta} (t))$ by \eqref{eq:LyaDrift}, and sum \eqref{ineq:cmp_f1} over $t = 1, 2, \dots, T$, yielding
\begin{equation} \label{ineq:cmp_f2}
    \begin{aligned}
        & \mathbb{E} \left[ L(\boldsymbol{\Theta}(T+1)) - L(\boldsymbol{\Theta} (1)) \right] + V \sum_{t = 1}^T \mathbb{E} \left[ \hat{f} (t) \right] \\
        \le ~& MT + V \sum_{t = 1}^T \mathbb{E} \left[ f^* (t) \right] - \sum_{g = 1}^G \sum_{t = 1}^T \mathbb{E} \left[ (q_g(t) + z_g(t)) x^*_g(t) \right] \\
        & + \frac{1}{2} \sum_{g = 1}^G \sum_{t = 1}^T \mathbb{E} \left[ (x_g^*(t))^2 \right]
    \end{aligned}
\end{equation}
According to \eqref{ineq:xgub}, $x_g^*(t)$ is bounded, then $(x_g^*(t))^2$ is also bounded by
\begin{equation} \label{ineq:xg_sq}
    0 \leq (x_g^*(t))^2 \leq (\overline{X}_g)^2, ~ \forall g, ~ \forall t.
\end{equation}
Divide both sides of \eqref{ineq:cmp_f2} by $VT$ and let $T$ goes to infinity, then we have
\begin{equation}
    \begin{aligned}
        \lim_{T \to \infty} \frac{1}{T} \sum_{t = 1}^T \mathbb{E} \left[ \hat{f} (t) - f^* (t) \right] 
        %\le ~ & \frac{M}{V} + \frac{1}{2V} \sum_{g = 1}^G \lim_{T \to \infty} \frac{1}{T} \sum_{t = 1}^T \mathbb{E} \left[ (x_g^*(t))^2 \right], \\
        \le \frac{M}{V} + \frac{1}{2V} \sum_{g = 1}^G (\overline{X}_g)^2,
    \end{aligned}
\end{equation}
where the right-hand side is reduced by \eqref{eq:mrs_q}, \eqref{eq:mrs_z} and \eqref{ineq:xg_sq}.
\end{proof}

\renewcommand{\theequation}{B.\arabic{equation}}

\section{Proof of Proposition \ref{prp:delay}} \label{appd:delay}

\begin{proof}    
Suppose Proposition \ref{prp:delay} does not hold, i.e., there exists charging demand $a_g(t_0)$ that cannot be served on or before the time slot $t_0 + D_g$. Since $q_g(t)$ and $z_g(t)$ are processed by FIFO method, we have the following conditions:
\begin{align}
    & q_g(t) > 0, ~ \forall t \in [t_0, t_0 + D_g], \label{ineq:qgpos} \\
    & \sum_{t = t_0}^{t_0 + D_g} x_g(t) < Q_g. \label{ineq:xq_Qg}
\end{align}
Then \eqref{eq:zgt} can be simplified by \eqref{ineq:qgpos}:
\begin{equation}
    z_g(t + 1) \ge z_g(t) + \frac{\alpha_g}{R_g} - x_g(t), ~ \forall t \in [t_0, t_0 + D_g]. \label{ineq:zgpos}
\end{equation}
Summing \eqref{ineq:zgpos} over $[t_0, t_0 + D_g]$, we obtain
\begin{equation}
    z_g(t_0 + D_g + 1) - z_g(t_0) \ge \frac{\alpha_g D_g}{R_g} - \sum_{t = t_0}^{t_0 + D_g} x_g(t). 
\end{equation}
Obviously, $z_g(t) \ge 0, ~ t = 1, 2, \dots, T$, thus we have
\begin{align}
    Z_g \ge z_g(t_0 + D_g + 1) - z_g(t_0) > \frac{\alpha_g D_g}{R_g} - Q_g,
\end{align}
where the second inequality is due to \eqref{ineq:xq_Qg}. Hence,
\begin{align}
    D_g < \frac{\alpha_g (Q_g + Z_g)}{R_g},
\end{align}
which contradicts the definition of $D_g$ \eqref{ineq:delay}. That completes the proof. 
\end{proof}

\renewcommand{\theequation}{C.\arabic{equation}}

\section{Proof of Proposition \ref{prp:con}} \label{appd:con}

\begin{proof}
For simplicity, the notation $t$ is omitted here since the proof can be applied to any time slot. % First, we will show that the optimal function $u(x)$ is continuous. 
Denote $x = r(\pi), ~ x_g = r_g(\pi), ~ g = 1, 2, \dots, G$, where $r_g(\pi)$ is a piecewise linear function given in \eqref{eq:xg1} or \eqref{eq:xg2}. According to \eqref{eq:x_def} and \eqref{eq:u_def}, we have
\begin{align}
    r(\pi) & = \sum_{g = 1}^G r_g(\pi), \\
    u(x) & = \int_0^{x} h(\xi) \text{d} \xi = \int_0^{x} r^{-1}(\xi) \text{d} \xi, \label{eq:uhr}
\end{align}
Obviously, $r_g(\pi)$ is continuous on $[0, +\infty)$.
% $h(x) \in (-\infty, +\infty), ~ \forall x \in \left[0, \sum_{g = 1}^G X_g \right]$, 
$u(x)$ is continuous since the following equality holds for any $x_0 \in \left[0, \sum_{g = 1}^G X_g \right]$:
\begin{equation}
    \lim_{x \to x_0^-} u(x) = \int_0^{x_0} h(x) \text{d} x = \lim_{x \to x_0^+} u(x).
\end{equation}
$r_g(\pi)$ is decreasing on $[0, +\infty)$ according to \eqref{eq:xg1} and \eqref{eq:xg2}, thus $h(x)$ is decreasing and $u(x)$ is concave on $\left[0, \sum_{g = 1}^G X_g \right]$ by \eqref{eq:uhr}.
\end{proof}

\renewcommand{\theequation}{D.\arabic{equation}}

\section{Proof of Proposition \ref{prp:agg}} \label{appd:agg}

\begin{proof}
For simplicity, the notation $t$ is omitted here since the proof can be applied to any time slot. By strong duality, solving \textbf{P3} is equivalent to solving the following problem:
\begin{align} \label{prbl:dual}
     \mathop{\max}_{\lambda \atop \underline{\chi}_l, \overline{\chi}_l \in \mathbb{R}^+} \mathop{\min}_{\underline{P}_i \le p_i \le \overline{P}_i \atop 0 \leq x_k \leq X_k} \mathcal{L}(\boldsymbol{p},\boldsymbol{x})
\end{align}
where $\boldsymbol{p} = (p_1, p_2, \dots, p_I), ~ \boldsymbol{x} = (x_1, x_2, \dots, x_K)$, $\mathcal{L}(\boldsymbol{p},\boldsymbol{x})$ is the Lagrangian function associated with \eqref{prbl:dual} given by
\begin{equation}
\begin{aligned}
    \mathcal{L}(\boldsymbol{p},\boldsymbol{x}) = & \sum_{i = 1}^I f_i^{\text{Gen}}(p_i) - \sum_{k = 1}^K u_k(x_k) \\
    & + \lambda \left(\sum_{j = 1}^J d_j + \sum_{k = 1}^K x_k - \sum_{i = 1}^I p_i \right) \\
    & + \sum_l \left( \underline{\chi}_l (-F_l - p_l^{\text{line}}) + \overline{\chi}_l (p_l^{\text{line}} - F_l) \right),
\end{aligned}
\end{equation}
$p_l^{\text{line}}$ is given in \eqref{eq:p_line}.
% Then the electricity price taken by the EV aggregator $k$ is given by
% \begin{equation} \label{eq:price1}
%     \pi_k = \frac{\partial \mathcal{L}}{\partial x_k} 
%     = \lambda + \sum_l \left( \delta_{kl} \underline{\chi}_l - \delta_{kl} \overline{\chi}_l \right), ~ \forall k,
% \end{equation}
\eqref{prbl:dual} can be decomposed into $I + K$ sub-problems solved by generators and aggregators independently. Specifically, the aggregator $k$ solves the following problem:
\begin{equation} \label{prbl:dual_k}
    \mathop{\max}_{\lambda \atop \underline{\chi}_l, \overline{\chi}_l \in \mathbb{R}^+} \min_{0 \leq x_k \leq X_k} -u_k(x_k) + \lambda x_k + x_k \sum_l \left( \delta_{kl} \underline{\chi}_l - \delta_{kl} \overline{\chi}_l \right)
\end{equation}
For simplicity, the subscript $k$ will be omitted in the following discussion since the following proof can be applied to any other aggregator.
Taking \eqref{eq:price} into \eqref{prbl:dual_k}, the inner-level problem is reformulated as:
\begin{equation} \label{prbl:dual_k_new}
    \min_{0 \leq x \leq X} -u(x) + \pi x
\end{equation}
Let $\hat{x}$ be the optimal charging power of all groups obtained by solving \textbf{P2}, i.e., $\hat{x} = \sum_{g = 1}^G \hat{x}_g$. Obviously, $0 \leq \hat{x} \leq X$. Moreover, we have the following equation based on \eqref{eq:h_def} and \eqref{eq:u_def}:
\begin{equation}
    -\frac{\text{d} u}{\text{d} x} \bigg|_{x = \hat{x}} + \pi = 0,
\end{equation}
which is exactly the optimality condition of \eqref{prbl:dual_k_new}. Hence, $\hat{x}$ is the optimal solution of \eqref{prbl:dual_k_new}.

% Next, we discuss the optimality conditions of \eqref{prbl:dual_k_new} case by case.
% \begin{enumerate}
%     \item If $\hat{x}_k = 0$
%     \item If $\hat{x}_k = X_k$
%     \item If $\hat{x}_k \in (0, X_k)$, thus $\hat{x}_k$ satisfies
%     \begin{equation} \label{eq:xk_hat}
%         \frac{\text{d} u_k}{\text{d} x_k} \bigg|_{x_k = \hat{x}_k} = \pi_k.
%     \end{equation}
%     Let $u_k(x) = \int_0^x \pi_k \text{d}x$, then \eqref{eq:xk_hat} is guaranteed.
%     Recall that $h(x)$ in \eqref{eq:u_def} is derived by \eqref{eq:xg1} or \eqref{eq:xg2}, which is the optimal solution of \textbf{P2}. Hence, the optimal solution of \textbf{P2} and \textbf{P3} are equivalent.
% \end{enumerate}

% where $\lambda(t)$, $\underline{\chi}_l(t)$ and $\overline{\chi}_l(t)$ are the dual variables corresponding to the constraints \eqref{eq:ErgBal}, the lower and upper bounds of $p_l^{\text{line}} (t)$ in \eqref{ineq:LinePwr} respectively.
\end{proof}

\renewcommand{\theequation}{E.\arabic{equation}}

\section{Linear Counterpart of the Proposed Algorithm} \label{appd:lin}

Note that the upper bound of the Lyapunov drift \eqref{ineq:dpp} is different from standard Lyapunov optimization framework in the literature since the former contains the quadratic term $x_g^2(t)$. Intuitively, the modification may offer a range of benefits, such as improved optimality, except for increasing computational burdens.
To compare the performance of these two methods, we derive the linear counterpart of \textbf{P2} in this section.

\begin{figure}[ht]
  \centering
  \includegraphics[width=0.45\textwidth]{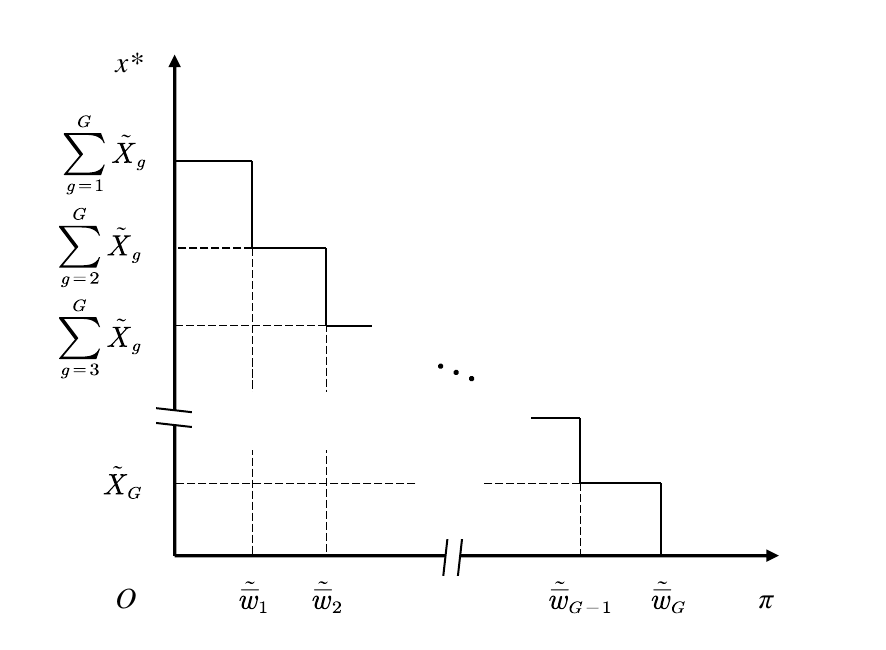} \\
  \caption{Bidding curve of a price-taking EV aggregator.} \label{fig:bidding_curve}
\end{figure}

By replacing $x_g^2(t)$ in \eqref{ineq:Q2-Ub0} by its upper bound $X_g^2$, we have a simpler form of \eqref{ineq:Q2-Ub}:
\begin{equation} \label{ineq:Q2-Ub1}
\begin{aligned}
& \frac{1}{2} \left(q_g^2(t+1) - q_g^2(t)\right) \\
\le ~& \frac{1}{2} \left(A_g^2 + X_g^2\right) 
+ q_g(t) \left(a_g(t) - x_g(t)\right)
\end{aligned}
\end{equation}
Finally, we get rid of the quadratic term in the objective of \textbf{P2}, yielding the following linear programming (LP) problem:
\begin{align}
  \textbf{P2$'$:}~ \min_{x_g(t),\forall g} ~ & \sum_{g = 1}^G (V\pi(t)-q_g(t) - z_g(t)) x_g(t) \label{obj:EVLin} \\
  \hbox{s.t.}~ & \eqref{ineq:xg}. \nonumber
\end{align}
The solution of \textbf{P2}$'$ is given by:
\begin{align} \label{eq:xg_prime}
    x'_g(t) = \left\{\begin{array}{ll}
      X_g, & \text{if}~ \pi(t) \in \left[0, \overline{w}_g(t) \right],\\      
      0, & \text{if}~ \pi(t) \in \left( \overline{w}_g(t), +\infty \right).
    \end{array} \right.
\end{align}
We sort the elements of $(\overline{w}_g)_{g=1}^G$ in ascending order and denote the new sequence by $\left( \tilde{\overline{w}}_g \right)_{g=1}^{G}$, i.e., $\tilde{\overline{w}}_1 \le \tilde{\overline{w}}_2 \le \dots \le \tilde{\overline{w}}_{G}$. The maximum charging power of the group corresponding to $\tilde{\overline{w}}_g$ is denoted by $\tilde{X}_g$.
Then the total charging power is given by:
\begin{align} \label{eq:x_prime}
    x' = \left\{\begin{array}{ll}
      \sum_{g=1}^G \tilde{X}_g, & \text{if}~ \pi \in [0, \tilde{\overline{w}}_1],\\
      \dots & \\
      \sum_{g=n}^G \tilde{X}_g, & \text{if}~ \pi \in (\tilde{\overline{w}}_{n-1}, \tilde{\overline{w}}_n],\\
      \dots & \\
      0, & \text{if}~ \pi \in (\tilde{\overline{w}}_{G}, +\infty),
    \end{array} \right.
\end{align}
The optimal bidding curve is a stepwise curve shown in Fig. \ref{fig:bidding_curve}.

Hence, the optimal bidding cost curve is also a piecewise linear function, given by
\bsq
\begin{align} \label{eq:u_prime}
    u'(x) = \left\{\begin{array}{ll}
      \tilde{\overline{w}}_{G} x, & \text{if}~ x \in \left[0, \tilde{X}_G \right],\\
      u'_n(x), & \text{if } x \in \left( \sum\limits_{g=n}^G \tilde{X}_g, \sum\limits_{g=n-1}^G \tilde{X}_g \right], \\
      \dots & \\
      \sum\limits_{g=1}^G \tilde{\overline{w}}_{g} \tilde{X}_g, & \text{if}~ x \in \left( \sum\limits_{g=1}^G \tilde{X}_g, +\infty \right).
    \end{array} \right.
\end{align}
where
\begin{equation} \label{eq:u_n_prime}
    u'_n(x) = \sum_{g=n}^G \tilde{\overline{w}}_{g} \tilde{X}_g + \tilde{\overline{w}}_{n-1} \left(x - \sum_{g=n}^G \tilde{X}_g \right), ~ n = 2, \dots, G.
\end{equation}
\esq
% \begin{align} \label{eq:u_prime}
%     u(x) = \left\{\begin{array}{ll}
%       \tilde{\overline{w}}_{G} x, & \text{if}~ x \in \left[0, \tilde{X}_G \right],\\
%       \tilde{\overline{w}}_{G} \tilde{X}_g & \multirow{2}{*}{if $x \in \left( \tilde{X}_{G}, \sum\limits_{g = G-1}^G \tilde{X}_g \right]$,} \\
%       + \tilde{\overline{w}}_{G-1} \left(x - \tilde{X}_G \right), & \\
%       \sum\limits_{g=n}^G \tilde{\overline{w}}_{g} \tilde{X}_g & \multirow{2}{*}{if $x \in \left( \sum\limits_{g=n}^G \tilde{X}_g, \sum\limits_{g=n-1}^G \tilde{X}_g \right]$,} \\
%       + \tilde{\overline{w}}_{n-1} \left(x - \sum\limits_{g=n}^G \tilde{X}_g \right), & \\
%       \dots & \\
%       \sum\limits_{g=1}^G \tilde{\overline{w}}_{g} \tilde{X}_g, & \text{if}~ x \in \left( \sum\limits_{g=n}^G \tilde{X}_g, +\infty \right).
%     \end{array} \right.
% \end{align}
Specifically, there will be no $u_n(x)$ in \eqref{eq:u_prime} if $G = 1$. Note that $u(x)$ is actually the minimum of a group of affine functions in this case. Therefore, \textbf{P3} can be replaced by an LP problem as follows by taking \eqref{eq:u_prime} into \eqref{obj:OPF}:
\bsq
\begin{align}
    \textbf{P3$'$:} ~ \min_{p_i, x_k, \sigma_k}~ & \sum_{i = 1}^I f_i^{\text{Gen}}(p_i) - \sum_{k = 1}^K \sigma_k \\
    \mbox{s.t.}~ & \sigma_k \le \tilde{\overline{w}}_{k G_k} x_k, ~ \forall k, \\
    & \sigma_k \le u'_{kn}(x_k), ~ n = 2, \dots, G_k, \forall k, \\
    & \sigma_k \le \sum_{g=1}^{G_k} \tilde{X}_{kg}, ~ \forall k, \\
    & \eqref{eq:ErgBal} - \eqref{ineq:LinePwr}, \nonumber
\end{align}
\esq
where $\left( \tilde{\overline{w}}_{kg} \right)_{g=1}^{G_k}$ and $\tilde{X}_{kg}$ are the sorted sequence and the corresponding maximum charging power of the aggregator $k$ defined earlier, $u'_{kn}(x_k)$ is $u'_n(x)$ of the aggregator $k$ defined in \eqref{eq:u_n_prime}.

\end{document}